\documentclass[11pt,final]{article}

\usepackage{amsfonts,amsthm,amsmath,amssymb}
\usepackage{graphicx}
\usepackage{color}
\usepackage{hyperref}
\usepackage{a4wide}
\usepackage{multirow}
\usepackage{enumitem}
\usepackage[normalem]{ulem} 

\newtheorem{proposition}{Proposition}[section]

\newtheorem{lemma}[proposition]{Lemma}

\numberwithin{equation}{section}
\numberwithin{figure}{section}
\numberwithin{table}{section}

\allowdisplaybreaks[4]

\renewcommand{\phi}{\varphi}
\renewcommand{\theta}{\vartheta}
\renewcommand{\epsilon}{{\color{blue}\varepsilon}}

\renewcommand{\d}[1]{\mathinner{\mathrm{d}{#1}}}

\newcommand{\reali}{{\mathbb{R}}}

\newcommand{\caratt}[1]{\chi_{\strut#1}}

\begin{document}

\title{Infectious Disease Spreading Fought by \\
  Multiple Vaccines Having a Prescribed Time Effect}

\author{Rinaldo M.~Colombo$^1$ \and Mauro
  Garavello$^{2,}$\footnote{Corresponding author}}

\maketitle

\footnotetext[1]{INdAM Unit \& Department of Information Engineering,
  University of Brescia, Italy.}

\footnotetext[2]{Dept.~of Mathematics and its Applications, University
  of Milano - Bicocca}

\begin{abstract}
  We propose a framework for the description of the effects of
  vaccinations on the spreading of an epidemic disease.  Different
  vaccines can be dosed, each providing different immunization times
  and immunization levels. Differences due to individuals' ages are
  accounted for through the introduction of either a continuous age
  structure or a discrete set of age classes. Extensions to gender
  differences or to distinguish fragile individuals can also be
  considered. Within this setting, vaccination strategies can be
  simulated, tested and compared, as is explicitly described through
  numerical integrations.

  \medskip

\end{abstract}

\section{Introduction}
\label{sec:I}

We propose a modeling framework to simulate the global process of a
vaccination campaign to fight the spreading of an epidemic. Vaccines,
possibly with different characteristics, are dosed to susceptible
individuals. Each vaccine is identified by the efficiency and the
duration of the protection it provides. In our model also individuals
that recovered from the disease are immunized for a prescribed time
period, after which they get back to be susceptible. In the age
structured version, these times and efficiencies are assumed to be age
dependent.

A common strategy to insert vaccination and in particular the loss of
immunization in a SIR type model consists in assigning to these
phenomena a \emph{rate}, typically proportional to the number of
susceptible and vaccinated individuals. On the contrary, here we
introduce the time at which individuals are vaccinated, account for
the time dependent level of immunization provided by the vaccine and
give a precise role to the duration of this immunization.

The proposed class of models relies on a deterministic and macroscopic
description, developed on top of the SIR model, and displays an
evolution which is inherently \emph{``multiscale''}: a first time
scale is that of the pathogen diffusion, which interacts at different
time scales with the different vaccines and with the recovering from
the disease. For a stochastic approach, we refer for instance also
to~\cite{BertagliaLiuParechi}, while a fuzzy approach is
in~\cite{Al-QanessEtAl, RegisEtAl}. On the basis of the epidemic
evolution described by the present model, consequences at the social
or economic levels can be described as in~\cite{AlbiBertagliaEtAl,
  FabbriGozziZanco}, for instance, or~\cite{MR4385929, MR4147945}
where a kinetic model of wealth exchange is proposed. A summary of the
historical development of macroscopic models for virus diffusion and
vaccination is in~\cite{MR3821682}.

The interaction among the different populations, e.g., susceptible,
infected, vaccinated and recovered, combined with the different time
scales leads to the formation of oscillations or \emph{epidemic
  waves}~\cite{LemonMahmoud2005}. When no vaccination is dosed (e.g.,
Figure~\ref{fig:no0Vacc}), or even more when a very heavy vaccination
campaign is in place (e.g., Figure~\ref{fig:q1}), then these waves
fade out rather quickly. On the contrary, a relatively mild
vaccination campaign hinders the virus propagation without stopping
it, so that these epidemic waves become rather persistent (e.g.,
Figure~\ref{fig:q2}).

The present model allows to test/compare different vaccination
strategies. For instance, analyzing the number of casualties resulting
from a vaccination campaign that leaves a fixed percentage, say $S_*$,
of non vaccinated individuals shows a sort of \emph{``herd
  immunity''}~\cite{Randolph2020737} effect. Indeed, the number of
casualties suffers a sharp increment in correspondence to a threshold
value $S_*$, roughly close to $10\%$ of the initial population (see
Figure~\ref{fig:NoVax}).

The choice of the vaccination strategy gets even more relevant when
different vaccines are available. It is realistic to imagine that
different vaccines provide different levels of immunity for different
time periods~\cite{Kai2021221, Mukhopadhyay202193}. Then, for
instance, the use of a poor vaccine has a doubly negative
effect. First, it does not ensure a good level of immunization and,
second, may prevent vaccinated individuals to get a better vaccine as
long as its effect is in place, see~\S~\ref{subsec:comp-use-diff}.

Age differences, too, require careful planning of vaccination
campaigns. Consider for simplicity $2$ classes: \emph{``younger''}
individuals are more infective, while \emph{``older''} ones are more
fragile. A vaccination strategy consisting in dosing exclusively the
older ones first is not necessarily the best choice. Indeed, a
campaign where the proportions of young and old dosed is carefully
chosen according to the disease diffusion can reduce the number of
casualties, even in the old class, see
\S~\ref{subsec:age-depend-vacc}.

In the realizations of the present framework discussed below, we keep
on purpose the number of populations to a minimum. It goes without
saying that the extension to richer structures is easily achievable at
the cost of only technical complications. The current literature
provides various examples of multispecies/multicompartment models,
often compared with real measurements, see for
instance~\cite{Giordano2020855, MR4334820, Yang20202708}.

We stress that the setting here introduced is amenable to consider,
for instance, also movements in space, gender differences or the
presence of more fragile individuals. These extensions, clearly,
formally complicate the equations. However, their numerical treatment
fits in the brief description in Appendix~\ref{sec:append-note-numer}
and does not require the introduction of new or \emph{ad hoc}
algorithms. Movements in space can be comprised with a procedure
similar to that used in Section~\ref{sec:Gen} to introduce a
continuous age structure, possibly introducing a further distinction
among individuals having different destinations, see~\cite{preprint,
  FrancescaElena} for further details. A different approach to
diffusion in space is treated, for instance, in~\cite{MR3881874}. The
setting therein is based on stochastic ordinary differential
equations, eventually leading to partial differential equations of
second order in the space
derivative~\cite[Formula~(13)]{MR3881874}. For a discussion of gender
and age differences see~\cite[Table~1]{Russo20212517}.

Aiming at a quantitative fitting with specific data reasonably
requires to let the various functions and parameters defining the
evolution (e.g., recovery rate, vaccine's efficiency or duration,
infectivity, $\ldots$) depend on time. The introduction of time
dependencies may account, for instance, for seasonal effects, changes
in lockdown policies, improvements in drug efficacy, $\ldots$
see~\cite{MR4009539, Merow202027456}.

\medskip

The next section is devoted to the simplest case of a single vaccine
(whose effect has a prescribed duration) without age structure. Then,
Section~\ref{sec:Two} deals with the concurrent use of multiple
vaccines. Age structures, both continuous and discrete, are the
subject of Section~\ref{sec:Gen}: here, in particular, the effects of
vaccines depend on age. In Section~\ref{sec:real-data-fitting} we
address the issue of choosing proper values for the parameters and
functions in the models we introduced, on the basis of Covid--19
related data, mostly related to the Italian situation, from the
current literature.

\section{A Single Vaccine}
\label{sec:One}

We present here our framework in its simplest realization, namely
considering a single vaccine and we test different vaccination
strategies to control the spreading of the disease.

As a starting point~\cite[Formula~(5)]{MR3821682}, consider the SIR
model
\begin{equation}
  \label{eq:4}
  \left\{
    \begin{array}{l}
      \dot S = -\rho_S \, I \, S
      \\
      \dot I = \rho_S \, I \, S - \theta\, I - \mu \, I
      \\
      \dot R = \theta \, I
    \end{array}
  \right.
\end{equation}
where, as usual $S,I,R$ are the number (or percentages) of Susceptible
($S$), Infected ($I$) and Recovered ($R$) individuals. The infectivity
coefficient $\rho_S$, the recovery rate $\theta$ and the mortality
rate $\mu$ are here considered to be constant; were they time
dependent, only technical difficulties would arise. As is well known,
in~\eqref{eq:4} the total number of individuals varies, actually
diminishes, exclusively due to the mortality term, i.e.,
$\frac{\d{~}}{\d{t}} (S+I+R) = -\mu \, I$. When long time intervals
are considered, it might be appropriate to include mortality also in
the $S$ and $R$ equations, or also natality, typically only in the $S$
equation. Other realizations might comprehend also time dependent
immigration/emigration terms, for instance.

As a first step, we modify~\eqref{eq:4} to allow for recovered
individuals to get re-infected, after a time $T_R$ from recovery. To
this aim, we modify the unknown $R$ to $R = R (t,\tau)$, the variable
$\tau$ being the time since recovery, with $\tau \in [0, T_R]$:
\begin{equation}
  \label{eq:5}
  \left\{
    \begin{array}{l}
      \dot S = -\rho_S \, I \, S + R (t,T_R)
      \\
      \dot I = \rho_S \, I \, S + \int_0^{T_R} \rho_R (\tau) \, R (t,\tau) \d\tau \, I - \theta\, I - \mu \, I
      \\
      \partial_t R + \partial_\tau R = - \rho_R \, R \, I
      \\
      R (t,0) = \theta \, I \,.
    \end{array}
  \right.
\end{equation}
Here, the $R$ compartment displays an \emph{``internal dynamics''},
see~\cite{FrancescaElena}. In other words, $R (t,\tau)$ is the number
of individuals at time $t$ that recovered at time
$t-\tau$. Elementary, though useful, is to note that the $R$
in~\eqref{eq:5} and the variable bearing the same name in~\eqref{eq:4}
have different dimensions. As above, the total number of individuals
varies, namely diminishes, exclusively due to mortality, i.e.,
\begin{displaymath}
  \frac{\d{~}}{\d{t}} \left(S (t) + I (t) +\int_0^{T_R} R (t,\tau)\d\tau\right)
  =
  -\mu \, I (t) \,.
\end{displaymath}
The function $\rho_R = \rho_R (\tau)$ describes how easy/difficult it
is that an $R$ individual gets infected after time $\tau$ from
recovery. A possible reasonable behavior of the map
\begin{figure}[!h]
  \centering
  \begin{minipage}{0.5\linewidth}
    \includegraphics[width = \linewidth]{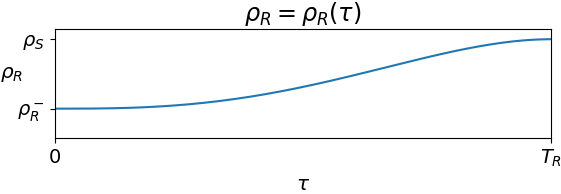}
  \end{minipage}%
  \begin{minipage}{0.5\linewidth}
    \begin{equation}
      \label{eq:7}
      \begin{array}{r@{\,}c@{\,}l@{}}
        \rho_R (\tau)
        & =
        & \rho_R^- + (\rho_S-\rho_R^-) \; \Phi \left(\frac{\tau}{T_R}\right)
        \\
        \Phi (s)
        & =
        & (4-3s) \, s^3
      \end{array}
    \end{equation}
  \end{minipage}
  \caption{\small Left, a reasonable choice of the function $\rho_R$:
    at $\tau\approx0$ we have $\rho_R (\tau) \approx \rho_R^-$, a
    \emph{small} ($\rho_R^- \ll \rho_S$) value quantifying the
    immunization resulting from recovering. As time from recovery
    passes, $\rho_R (\tau)$ increases and at time $T_R$ attains the
    value $\rho_S$ of susceptibles. Right, the actual expression used
    in the diagram on the left and in the numerical integrations in
    Section~\ref{sec:One}. The relevant properties of $\Phi$ are its
    continuity and monotonicity, from $\Phi (0) = 0$ to
    $\Phi (1)= 1$.\label{fig:rhoR}}
\end{figure}
$\tau \mapsto \rho_R (\tau)$ is depicted in Figure~\ref{fig:rhoR}.
For $\tau$ near to $0$, $\rho_R (\tau)$ equals $\rho_R^-$, a value far
smaller than the infectivity coefficient $\rho_S$ in~\eqref{eq:4}
or~\eqref{eq:5}. As the time $\tau$ from recovery grows, also
$\rho_R (\tau)$ grows and gets back to the value $\rho_S$ at time
$T_R$, when recovered individuals return to be susceptible. The
extension to $\rho_R$ depending also on $t$ is immediate, as also that
of letting $T_R \to +\infty$, as explicitly considered below.

In~\eqref{eq:5}, the rate $\theta\, I$ at which infected individuals
recover, tuned through the constant $\theta$, is the same as
in~\eqref{eq:4}. Each recovered individual after time $\tau = T_R$
from recovery gets back to being susceptible.

Remark that when considering a finite number of age classes
$a_1, a_2, \ldots, a_k$ or a continuous age structure with
$a \in \reali_+$, then different ages may well have different times
$T_R$, i.e., $T_R = T_R (a)$.

\smallskip

The effect of a vaccination that does not ensure permanent immunity is
to some extent similar to the temporary immunization of recovered
individuals as described above. A first difference is that
immunization is obtained some time after being dosed. More relevant,
vaccinations depend on a vaccination strategy, i.e., on the arbitrary
choice of which and how many susceptibles are dosed at each
time. Therefore, we introduce a new population, namely $V$, where
$V (t,\tau)$ is the number of vaccinated individuals at time $t$ that
were dosed at time $t-\tau$, so that $\tau$ here is the time since
vaccination.

\begin{figure}[!h]
  \centering \includegraphics[width = 0.75\linewidth]{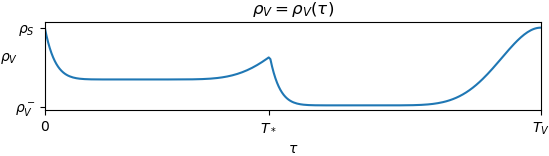}
  \caption{\small Qualitative behavior of a possible map
    $\tau \mapsto \rho_V (\tau)$ in the case of a vaccination
    consisting of $2$ doses, whose effect ceases at time $T_V$. The
    second shot takes place at time $T_*$ after the first one. After
    time $T_V$ from the first dose, the protection provided by the
    vaccine expires, since $\rho_V (T_V)$ attains the value
    $\rho_S$. It is not required that $\rho_V$ vanishes on an
    interval: the efficiency of the vaccine translates into $\rho_V$
    being \emph{''very small''} for a \emph{``very long''}
    time.\label{fig:2doses}}
\end{figure}

We are thus lead to introduce the model
\begin{equation}
  \label{eq:2}
  \left\{
    \begin{array}{l@{}}
      \displaystyle
      \dot S
      =
      - \rho_S \, I \, S
      + V (t,T_V)
      + R (t, T_R)
      - p (t, S, V, I, R)
      \\
      \displaystyle
      \partial_t V + \partial_\tau V
      =
      - \rho_V \, V \, I
      \hfill \tau \in [0, T_V]
      \\
      \displaystyle
      \dot I
      =
      \left(
      \rho_S \, S
      + \int_0^{T_V} \!\! \rho_V (\tau) \, V (t,\tau) \d\tau
      + \int_0^{T_R}  \!\!\rho_R (\tau) \, R (t,\tau) \d\tau
      - (\theta + \mu)
      \right) I \qquad
      \\
      \displaystyle
      \partial_t R + \partial_\tau R
      =
      - \rho_R \, R \, I
      \hfill \tau \in [0, T_R]
      \\
      \displaystyle
      V (t, 0) = p(t,S, V, I, R)
      \\
      \displaystyle
      R (t, 0) = \theta \, I \,,
    \end{array}
  \right.
\end{equation}
where $T_V$ is the time when the immunization provided by vaccination
terminates.  Similarly to what is described above with reference to
the function $\rho_R = \rho_R (\tau)$, now the function
$\rho_V = \rho_V (\tau)$ describes how easy/difficult it is for an
individual dosed at time $t-\tau$ to get infected at time $t$,
i.e.~after time $\tau$ from vaccination, for $\tau \in [0, T_V]$.
\begin{figure}[!h]
  \centering
  \begin{minipage}{0.5\linewidth}
    \includegraphics[width = \linewidth]{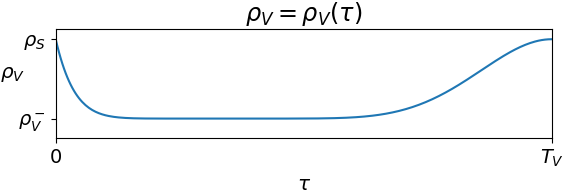}
  \end{minipage}%
  \begin{minipage}{0.5\linewidth}
    \begin{equation}
      \label{eq:8}
      \begin{array}{r@{\,}c@{\,}l@{}}
        \rho_V (\tau)
        & =
        & \rho_V^- + (\rho_S-\rho_V^-) \; \Psi \left(\frac{\tau}{T_V}\right)
        \\
        \Psi (s)
        & =
        & \left(1 - \frac{27}{4}s (1-s)^2\right)^4
      \end{array}
    \end{equation}
  \end{minipage}
  \caption{\small Left, a reasonable choice of the function $\rho_V$:
    at $\tau\approx0$ we have $\rho_V (\tau) \approx \rho_S$, since
    immunization is not immediate after being dosed. As time from
    vaccination passes, $\rho_V (\tau)$ decreases, reaches a lowest
    level $\rho_V^-$ and at time $T_V$ is back at the value $\rho_S$
    of susceptibles. Right, the actual expression used in the diagram
    on the left and in the numerical integrations in
    Section~\ref{sec:One}. The relevant properties of $\Psi$ are its
    continuity, the fact that $\Psi (0) = \Psi (1)= 1$ and the kind of
    \emph{plateau} near its minimum.\label{fig:rhoV}}
\end{figure}
Qualitatively, in the case of a vaccine consisting of a single dose,
the function $\rho_V$ can be chosen, for instance, as depicted in
Figure~\ref{fig:rhoV}. In the case of a vaccine consisting of $2$
shots, a possible behavior of $\tau \mapsto \rho_V (\tau)$ is in
Figure~\ref{fig:2doses}.
Quantitatively, in both cases, the function $\rho_V$ depends on
parameters specific to the vaccine under consideration.

In~\eqref{eq:2}, a key role is played by the function
$p = p (t, S, V, I, R)$. It describes the vaccination strategy,
quantifying how many susceptible individuals are dosed at time
$t$. Analytically, remark that the dependence of $p$ on the variables
$S, V, I, R$ may well be of a \emph{functional} nature, in the sense
that $p$ may depend, for instance, on time integrals of the functions
$S, V, I, R$, see~\eqref{eq:36}.

Several statistics on the solutions to~\eqref{eq:2} are of
interest. First, the total number of casualties $\mathcal{D} (t_0,T)$
between time $t_0$ and time $T$ (with $t_0 < T$) clearly equals the
variation in the total number of individuals between times $t_0$ and
$T$. It can be computed as
\begin{equation}
  \label{eq:27}
  \begin{array}{@{}r@{\,\,}c@{\,}l@{}}
    \mathcal{D} (t_0,T)
    & =
    & \displaystyle
      \left(
      S (t_0)
      + \int_0^{T_V} V (t_0,\tau)\d\tau
      + I (t_0)
      + \int_0^{T_R} R (t_0,\tau)\d\tau
      \right)
    \\
    &
    & \displaystyle
      -
      \left(
      S (T)
      + \int_0^{T_V} V (T,\tau)\d\tau
      + I (T)
      + \int_0^{T_R} R (T,\tau)\d\tau
      \right)
    \\
    & =
    & \displaystyle
      \int_{t_0}^T \mu \, I (t) \d{t}\,.
  \end{array}
\end{equation}
An estimate of the cost of the vaccination campaign is given by the
total number of vaccines dosed between time $t_0$ and time $T$, that
is
\begin{equation}
  \label{eq:35}
  \mathcal{V} (t_0,T)
  =
  \int_{t_0}^T p\left(t, S (t), V (t), I (t), R (t)\right) \d{t} \,.
\end{equation}
A common index used to measure the virus propagation is the basic
reproduction number~\cite[Section~10.2]{Murray1}, which is here
computed as
\begin{equation}
  \label{eq:3}
  \mathcal{R}_0 (t)
  =
  \dfrac{\rho_S \, S
    + \int_0^{T_V} \!\! \rho_V (\tau) \, V (t,\tau) \d\tau
    + \int_0^{T_R}  \!\!\rho_R (\tau) \, R (t,\tau) \d\tau}{\theta+\mu}
\end{equation}
since we have the equivalences
\begin{displaymath}
  \mathcal{R}_0 (t) > 1 \iff \dot I (t) > 0
  \quad \mbox{ and } \quad
  \mathcal{R}_0 (t) < 1 \iff \dot I (t) < 0 \,.
\end{displaymath}
Remark that the above expression of $\mathcal{R}_0 (t)$ does not
require the knowledge of the number of infected $I (t)$.

\subsection{Comparing Vaccination Strategies}
\label{subsec:comp1}

Our aim in the integrations below is to stress qualitative features of
the model~\eqref{eq:2}. Quantitative data are presented to allow the
reader to reproduce the results. Where helpful, we provide references
coherent with the quantitative choices adopted, bearing in mind that
several measurements are currently being improved and updated in the
literature. Nevertheless, it may help the reader to consider time as
measured in days, while $S (t)$, $I (t)$, $\int V (t,\tau) \d\tau$ and
$\int R (t,\tau) \d\tau$ are percentages, since the total initial
population is throughout fixed to $100$.

The numerical algorithm adopted is described in
Appendix~\ref{sec:append-note-numer}.

\begin{figure}[!h]
  \includegraphics[width = 0.33\linewidth]{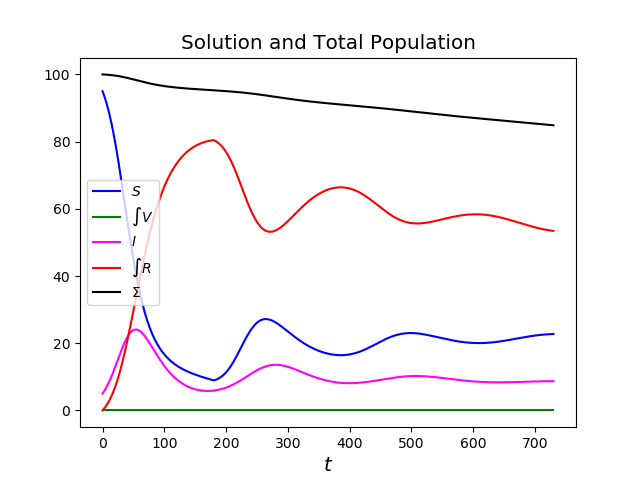}\hfil%
  \includegraphics[width = 0.33\linewidth]{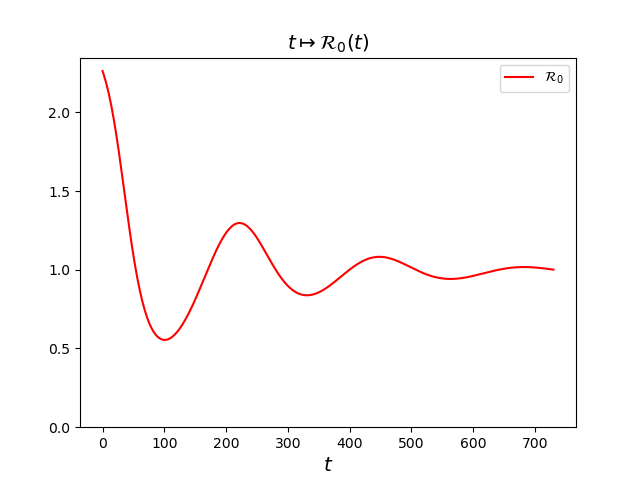}\hfil%
  \includegraphics[width = 0.33\linewidth]{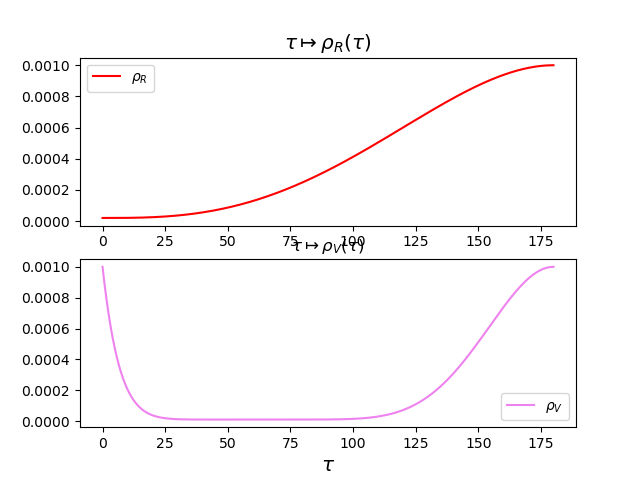}\\
  \caption{\small Solution
    to~\eqref{eq:2}--\eqref{eq:25}--\eqref{eq:26} with no
    vaccinations, i.e., $p\equiv 0$. Left, the graphs of
    $t \mapsto S (t)$, $t \mapsto \int_0^{T_V}V (t,\tau)\d\tau$,
    $t \mapsto I (t)$, $t \mapsto \int_0^{T_R}R (t,\tau)\d\tau$ and of
    their sum, labeled by $\Sigma$. Middle, the corresponding graph of
    $t \mapsto \mathcal{R}_0 (t)$, as defined in~\eqref{eq:3}. Right,
    the graphs of the functions $\tau \mapsto \rho_R (\tau)$, above,
    and of $\tau \mapsto \rho_V (\tau)$, below, used in this and in
    the forthcoming integrations. \label{fig:no0Vacc}}
\end{figure}

\paragraph{The Reference Situation:} We take as reference situation
the spreading of virus with no vaccination, described by~\eqref{eq:2}
with $p \equiv 0$ and with the following choices, which do not pretend
to be quantitatively fully justified by the available data:
\begin{equation}
  \label{eq:25}
  \begin{array}{r@{\,}c@{\,}l@{\quad\qquad}%
    r@{\,}c@{\,}l@{\quad\qquad}r@{\,}c@{\,}l}
    \rho_S
    & =
    & 1.0 \times 10^{-3}
    & \theta
    & =
    & 4.0 \times 10^{-2}
    & \mu
    & =
    & 2.0 \times 10^{-3}
    \\
    \rho_V^-
    & =
    & 1.0 \times 10^{-5}
    & \rho_V (\tau)
    & =
    & \rho_V^- + (\rho_S-\rho_V^-) \; \Psi (\tau/T_V)
    & T_V
    & =
    & 180
    \\
    \rho_R^-
    & =
    & 2.0 \times 10^{-5}
    & \rho_R (\tau)
    & =
    & \rho_R^- + (\rho_S-\rho_R^-) \; \Phi (\tau/T_R)
    & T_R
    & =
    & 180
  \end{array}
\end{equation}
where $\Psi$ is as in~\eqref{eq:8} and $\Phi$ is as
in~\eqref{eq:7}. Since $\theta/\mu = 20$, see~\cite{Russo20212517},
the above choice says that for an infected individual it is $20$ times
easier to recover than to die, corresponding to a mortality slightly
lower than $5\%$. The different partial immunizations provided by the
vaccine or by the recovery are described through the maps $\rho_V$ and
$\rho_R$, displayed in Figure~\ref{fig:no0Vacc}, right. In the
literature, available data keep being updated: with the present
choice~\eqref{eq:25}, the vaccine is more efficient than the
recovering, both because it leaves a lower probability to get infected
and because it is effective for a longer time.

The initial datum is
\begin{equation}
  \label{eq:26}
  S (0) = 95
  \,,\quad
  V (0,\tau) = 0 \mbox{ for } \tau \in [0, T_V]
  \,,\quad
  I (0) = 5
  \,,\quad
  R (0,\tau) = 0 \mbox{ for } \tau \in [0, T_R]
\end{equation}
meaning that at time $t=0$, the susceptibles are $95\%$ of the total
population, $5\%$ is infected, none is vaccinated and none is among
those who recovered.

In this reference situation, the casualties after time $730$ (i.e.,
$2$ years) amount to $15.1\%$ of the initial population. The numerical
integration shows the insurgence of \emph{``epidemic
  waves''}~\cite{LemonMahmoud2005}, see Figure~\ref{fig:no0Vacc},
left.

\bigskip

In the examples below, we always let the vaccination campaign begin
after time $t=30$, to allow for the onset of the virus spreading. This
is described through the term $\caratt{[30, +\infty[} (t)$ in the
vaccination strategy, see for instance~\eqref{eq:34}. Note also that
the time $T_V$ in~\eqref{eq:25} adopted below allows for multiple, up
to $4$, vaccinations of each single individual. Therefore, the number
of doses may well exceed the total initial population, set to $100$.

\paragraph{Leaving a Non Vaccinated Percentage:}
Practical considerations based on the different
attitudes~\cite{Wang20201} towards vaccines may induce or oblige to
avoid dosing a given portion of the population. Here, we describe this
situation through the vaccination strategy
\begin{equation}
  \label{eq:34}
  p (t, S, V, I, R)
  =
  \caratt{[30, +\infty[} (t) \; \caratt{]S_*,+\infty[} (S)
  \quad \mbox{ where }
  \begin{array}{r@{\,}c@{\,}l@{}}
    \caratt{[30, +\infty[} (t)
    & =
    & \left\{
      \begin{array}{l@{\;\mbox{ if }\;}r@{\,}c@{\,}l@{}}
        0
        & t
        & <
        & 30
        \\
        1
        & t
        & \geq
        & 30
      \end{array}
          \right.
    \\
    \caratt{]S_*, +\infty[} (S)
    & =
    & \left\{
      \begin{array}{l@{\;\mbox{ if }\;}r@{\,}c@{\,}l@{}}
        0
        & S
        & \leq
        & S_*
        \\
        1
        & S
        & >
        & S_*
      \end{array}
          \right.
  \end{array}
\end{equation}
meaning that when susceptibles are below the threshold value $S_*$,
the vaccination campaign stops. Note that, in the framework resulting
from~\eqref{eq:2}--\eqref{eq:34}, we do not impose that the non
vaccinated individuals are always the same.

\begin{figure}[!h]
  \centering \includegraphics[width =
  0.5\linewidth]{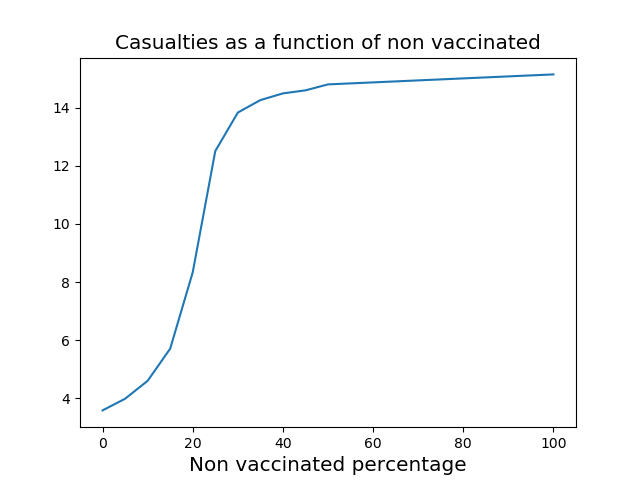}
  \caption{\small Casualties as a function of the non vaccinated
    percentage. Note the sharp increase starting already at a
    threshold of about $10\%$, somewhat justifying the term
    \emph{``herd immunity''}. The corresponding values are in
    Table~\ref{tab:novax}.}\label{fig:NoVax}
\end{figure}

As is to be expected, the higher the threshold $S_*$, the higher the
resulting number of casualties. However, we remark that when the
threshold percentage of non vaccinated gets near to $10\%$, the
corresponding number of casualties sharply increases, see
Figure~\ref{fig:NoVax}.
While it is somewhat arbitrary to choose a specific percentage where
this sharp increase begins, this behavior partly justifies the term
\emph{``herd immunity''}~\cite{Randolph2020737}, commonly used.

The actual computed values are in Table~\ref{tab:novax}, where the
case $S_* = 100$ corresponds to the reference situation above.
\begin{table}[!h]
  \centering
  \begin{tabular}{l|c|c|c|c|c|c|c|c|c|c|c|c}
    $S_*$
    & $0$
    & $5$
    & $10$
    & $15$
    & $20$
    & $25$
    & $30$
    & $35$
    & $40$
    & $45$
    & $50$
    & $100$
    \\ \hline
    Deaths
    & $3.59$
    & $3.99$
    & $4.60$
    & $5.71$
    & $8.34$
    & $12.5$
    & $13.8$
    & $14.3$
    & $14.5$
    & $14.6$
    & $14.8$
    & $15.1$
    \\ \hline
    Doses
    & 326
    & 298
    & 267
    & 227
    & 154
    & 61.3
    & 30.5
    & 17.7
    & 13.8
    & 11.3
    & 8.78
    & 0.00
  \end{tabular}
  \caption{\small Casualties and vaccinations in the solution
    to~\eqref{eq:2}--\eqref{eq:25}--\eqref{eq:26} corresponding to
    different values of $S_*$ in~\eqref{eq:34}. See also
    Figure~\ref{fig:NoVax}.\label{tab:novax}}
\end{table}

\paragraph{Automatic Feedback Based on $\mathcal{R}_0 (t)$:} Rather
than a systematic full speed vaccination campaign, as considered in
the preceding paragraph, one may consider a feedback strategy relying
on the index $\mathcal{R}_0$ defined in~\eqref{eq:3}. With the same
notation as in~\eqref{eq:34}, we set
\begin{equation}
  \label{eq:36}
  p (t,S, V, I, R)
  =
  p_* \;\;
  \caratt{[30, +\infty[} (t) \;\;
  \caratt{]0,+\infty[} (S) \;\;
  \caratt{]r_*,+\infty[} (\mathcal{R}_0)
\end{equation}
meaning that at time $t$, with $t > 30$, the campaign proceeds dosing
$p_*$ individuals per day, as soon as there are susceptibles (i.e.,
$S (t) > 0$) and $\mathcal{R}_0 (t)$ exceeds the threshold $r_*$.

This feedback strategy allows for a qualitative result, which is
\emph{independent} of the specific data and parameters chosen. Indeed,
assume the strategy $p$ in~\eqref{eq:36} is assigned so that
$\mathcal{R}_0$ is stabilized to $r_*$ after time $t_*$, i.e.,
$\mathcal{R}_0 (t) = r_*$ for $t \in [t_*, T]$ for a large $T$.  We
can clearly assume that $t_* > T_V$ and $t_* > T_R$. Then,
by~\eqref{eq:3}, the solution to model~\eqref{eq:2} for
$t \in [t_*, T]$ satisfies
\begin{equation}
  \label{eq:29}
  \left\{
    \begin{array}{rcl}
      \displaystyle
      I (t)
      & =
      & \displaystyle
        I (t_*) \; \exp\left((r_*-1) (\theta+\mu) (t-t_*)\right)
      \\
      \displaystyle
      R (t,\tau)
      & =
      & \displaystyle
        \theta \, I (t-\tau) \; \exp\left(-\int_{t-\tau}^t \rho_R \, I (s) \d{s}\right)
    \end{array}
  \right.
\end{equation}
see Lemma~\ref{lem:si}. As expected, in the case $r_* = 1$,
stabilizing $\mathcal{R}_0 (t)$ for $t \in [t_*,T]$, also $I$ is
stabilized at the value $I_* = I (t_*)$, and
$R (t,\tau) = \theta\, I_* \, e^{-\rho_R I_* \, \tau}$ is independent
of $t$. Note that casualties, defined in~\eqref{eq:27}, grow linearly
with time, proving that $T$ is necessarily bounded, its largest
possible value corresponding to when all individuals die.

\begin{figure}[!h]
  \includegraphics[width = 0.5\linewidth]{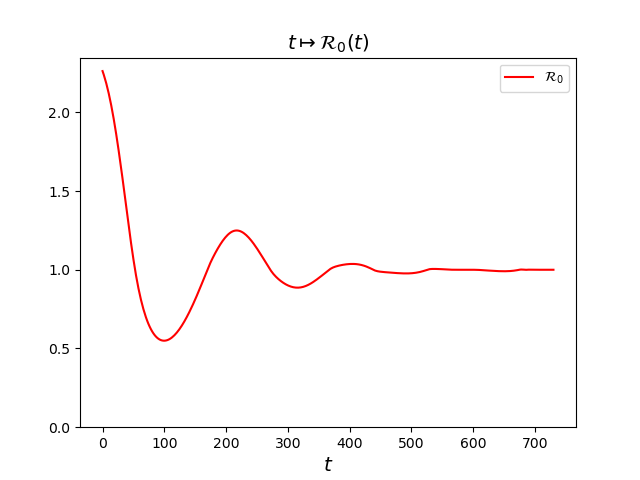}\hfil%
  \includegraphics[width = 0.5\linewidth]{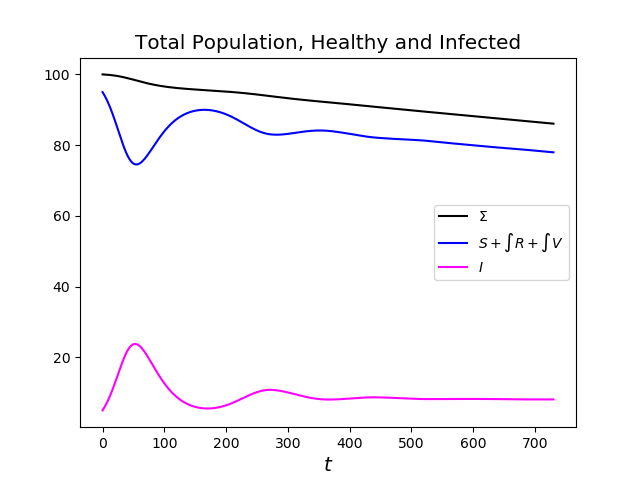}\hfil\\
  \includegraphics[width = 0.5\linewidth]{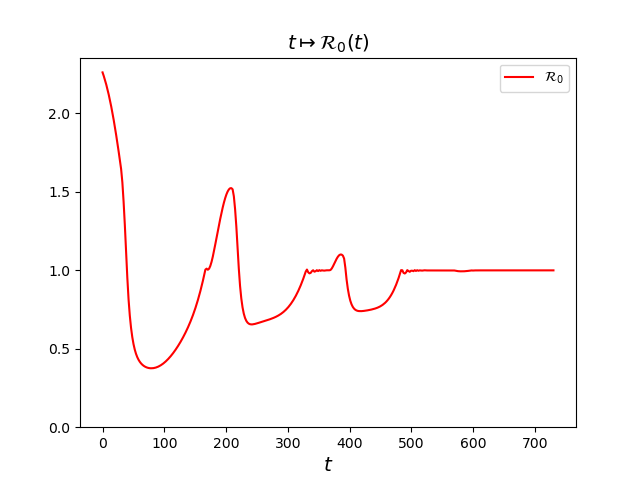}\hfil%
  \includegraphics[width = 0.5\linewidth]{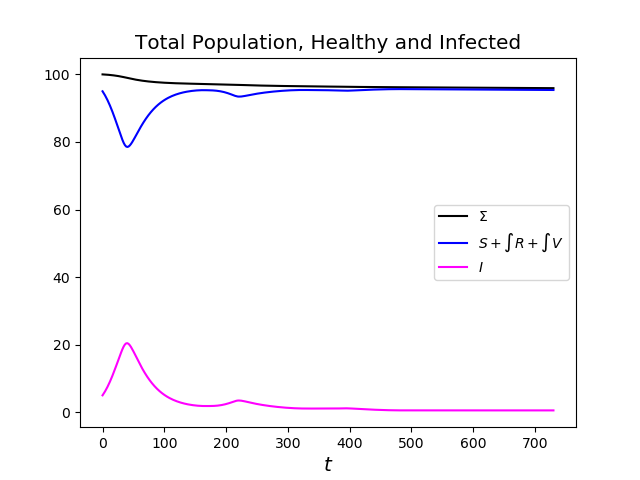}\\
  \caption{\small Solutions
    to~\eqref{eq:2}--\eqref{eq:25}--\eqref{eq:26} and the resulting
    values of $\mathcal{R}_0 (t)$ under strategy~\eqref{eq:36},
    $r_* =1.00$ with, above, $p_* = 0.1$ and, below, $p_* = 4.0$. The
    lower vaccination rate above causes a stabilization of $I$ at a
    higher value and, hence, a higher mortality.\label{fig:RRbis}}
\end{figure}

For arbitrary values of $r_*$, the former relation in~\eqref{eq:29}
immediately gives for $r_* \neq 1$,
\begin{equation}
  \label{eq:30}
  \mathcal{D} (t_*,T)
  =
  \dfrac{1}{1-r_*} \;
  \dfrac{\mu}{\theta+\mu}
  \left(
    1
    -
    e^{- (1-r_*) (\theta+\mu) (T-t_*)}
  \right)
  I (t_*) \,.
\end{equation}
Thus, for the disease to disappear, it is necessary to stabilize
$\mathcal{R}_0 (t)$ at a value $r_*$ strictly lower than $1$. However,
this condition is clearly not sufficient: one should also require that
$\mathcal{D} (t_*,T)$ in~\eqref{eq:30} does not exceed the number of
living individuals at time $t_*$.
\begin{figure}[!h]
  \includegraphics[width = 0.5\linewidth]{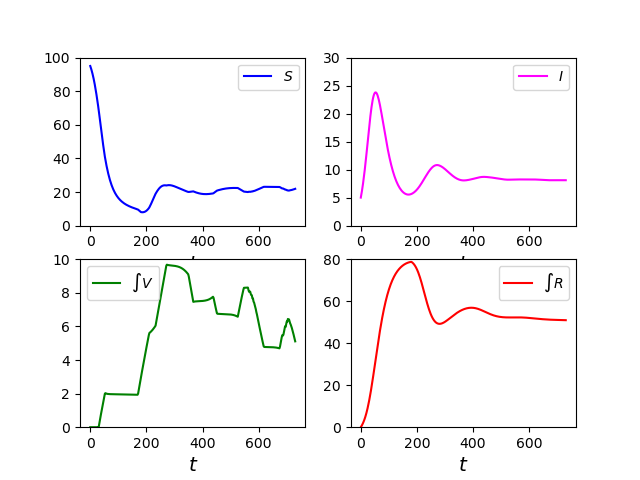}\hfil%
  \includegraphics[width = 0.5\linewidth]{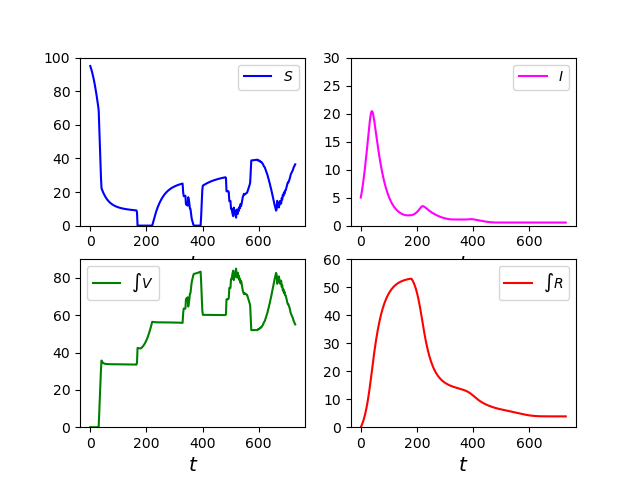}\hfil%
  \caption{\small Solutions
    to~\eqref{eq:2}--\eqref{eq:25}--\eqref{eq:26} with
    strategy~\eqref{eq:36} with $r_* = 1.00$ and, left, $p_* = 0.1$
    while, right, $p_* = 4.0$. Note that in both cases $I$ is
    stabilized at a strictly positive value, left $I_* \approx 7$ and,
    right, $I_* \approx 1.5$, causing casualties to keep increasing,
    see Table~\ref{tab:RR}. (The left and right scales in the lower
    diagrams of $\int V$ differ.)\label{fig:RR}}
\end{figure}
The integrations in Figure~\ref{fig:RR} confirm that stabilizing at
$\mathcal{R}_0 (t) = 1$ does not stop the spreading of the disease, as
also shown in Table~\ref{tab:RR}. When $\mathcal{R}_0 (t) = 1$, a sort
of \emph{``dynamic equilibrium''} is onset, so that, for large $t$,
the maps $t \to I (t)$ and $t \to \int R (t,\tau) \d\tau$ are
approximately constant, while $t \to S (t)$ and
$t \to \int V (t,\tau)\d\tau$ have oscillations that approximately
balance each other, so that their sum keep diminishing at a rate
approximately $\mu \, I (t_*)$, see Figure~\ref{fig:RRbis}.

\begin{figure}[!h]
  \includegraphics[width = 0.5\linewidth, trim=10 20 25 0]%
  {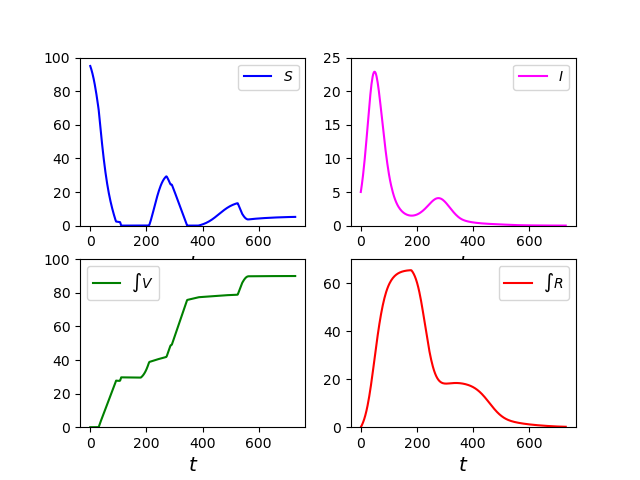}\hfil%
  \includegraphics[width = 0.5\linewidth, trim=10 20 25 0]%
  {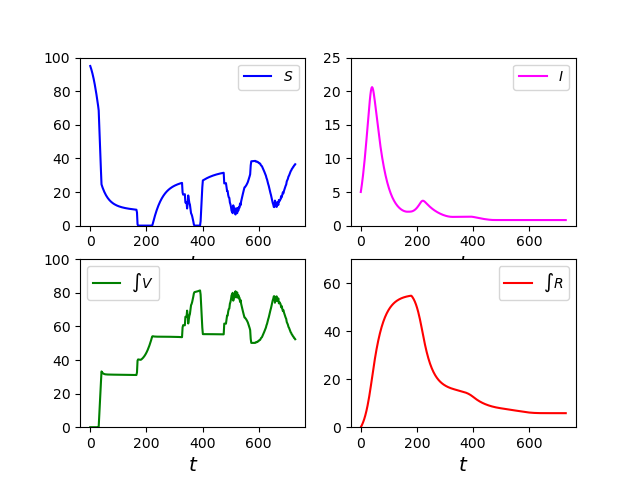}\hfil\\
  \includegraphics[width = 0.5\linewidth, trim=10 20 25 0]%
  {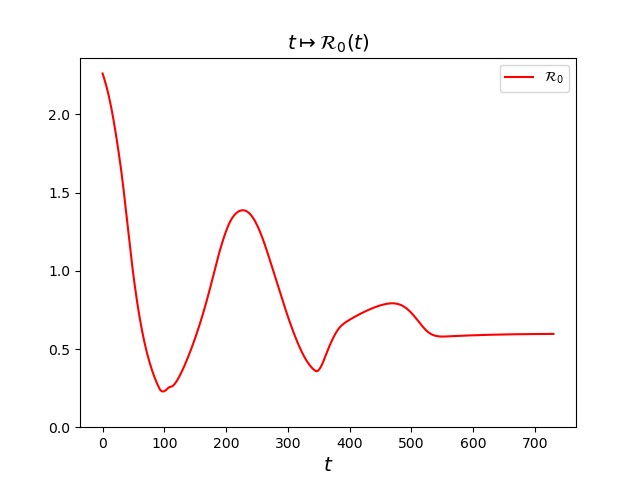}\hfil%
  \includegraphics[width = 0.5\linewidth, trim=10 20 25 0]%
  {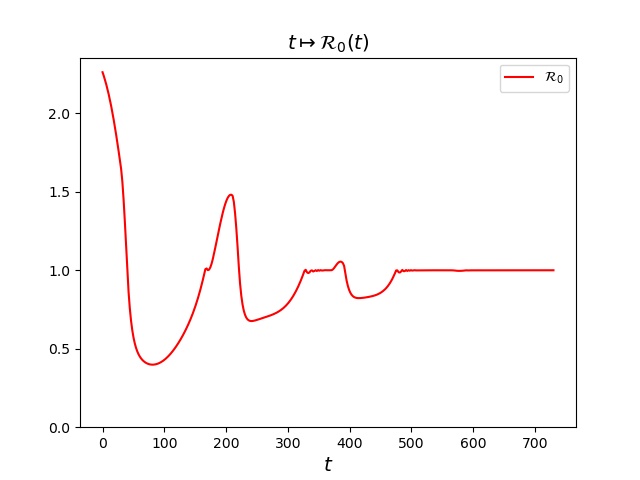}\\
  \caption{\small Solutions
    to~\eqref{eq:2}--\eqref{eq:25}--\eqref{eq:26} with
    strategy~\eqref{eq:36} and the resulting values of
    $\mathcal{R}_0 (t)$ corresponding to the bold data in
    Table~\ref{tab:RR}: left, $p_* = 0.50$ and $r_* = 0.25$; right,
    $p_*=3.50$ and $r_* = 1.00$; above, the solutions and, below, the
    value of $\mathcal{R}_0 (t)$ along the solutions. The higher
    vaccination rate in the $2$ diagrams on the right allows to get a
    lower number of casualties with a lower number of doses. In
    subsequent times, the lower value of $\mathcal{R}_0 (t)$ ensures
    that the choice on the left will be more effective in reducing
    the mortality.\label{fig:RRter}}
\end{figure}

Somewhat surprisingly, two situations arise where a higher vaccination
speed allows for a faster reduction of the infected and infectious
population, so that -- on the time interval considered -- the
resulting number of casualties is lower than that obtained after a
higher but slower vaccination campaign, see the bold data in
Table~\ref{tab:RR}.
\begin{table}[!h]
  \centering
  \begin{tabular}{cc|c|c|c}
    Deaths
    &
    & \multicolumn{3}{c}{Threshold $r_*$}
    \\
    $\mathcal{D} (0, 730)$
    &
    & 0.25
    & 0.50
    & 1.00
    \\ \hline
    \multirow{9}{*}{Speed $p_*$}
    & 0.10
    & 12.1
    & 12.1
    & 13.9
    \\ \cline{2-5}
    & 0.50
    & \textbf{4.65}
    & 4.89
    & 10.2
    \\ \cline{2-5}
    & 1.00
    & 3.63
    & 3.98
    & 8.61
    \\ \cline{2-5}
    & 1.50
    & 3.18
    & 3.53
    & 7.40
    \\ \cline{2-5}
    & 2.00
    & 2.94
    & 3.22
    & 6.44
    \\ \cline{2-5}
    & 2.50
    & 2.79
    & 2.99
    & 5.66
    \\ \cline{2-5}
    & 3.00
    & 2.68
    & 2.80
    & 5.07
    \\ \cline{2-5}
    & 3.50
    & 2.60
    & 2.64
    & \textbf{4.45}
    \\ \cline{2-5}
    & 4.00
    & 2.53
    & 2.55
    & 4.07
  \end{tabular}
  \qquad\qquad\qquad
  \begin{tabular}{cc|c|c|c}
    Doses
    &
    & \multicolumn{3}{c}{Threshold $r_*$}
    \\
    $\mathcal{V} (0,730)$
    &
    & 0.25
    & 0.50
    & 1,00
    \\ \hline
    \multirow{9}{*}{Speed $p_*$}
    & 0.10
    & 70.0
    & 69.0
    & 27.8
    \\ \cline{2-5}
    & 0.50
    & \textbf{291}
    & 287
    & 109
    \\ \cline{2-5}
    & 1.00
    & 325
    & 315
    & 145
    \\ \cline{2-5}
    & 1.50
    & 336
    & 328
    & 173
    \\ \cline{2-5}
    & 2.00
    & 341
    & 334
    & 195
    \\ \cline{2-5}
    & 2.50
    & 345
    & 338
    & 214
    \\ \cline{2-5}
    & 3.00
    & 348
    & 345
    & 229
    \\ \cline{2-5}
    & 3.50
    & 350
    & 349
    & \textbf{248}
    \\ \cline{2-5}
    & 4.00
    & 351
    & 351
    & 261
  \end{tabular}
  \caption{\small Casualties~\eqref{eq:27}, left, and total number of
    vaccinations~\eqref{eq:35}, right, corresponding to the
    strategy~\eqref{eq:36}
    in~\eqref{eq:2}--\eqref{eq:25}--\eqref{eq:26}. The bold data
    correspond to an unusual situation where less doses allow for less
    casualties, see also Figure~\ref{fig:RRter}.\label{tab:RR}}
\end{table}
More precisely, a higher vaccination speed allows for a faster
reduction of the $I$ population and to quickly dose all the $S$
individuals or lower $\mathcal{R}_0$ below the desired threshold, see
Figure~\ref{fig:RRter}.
On the other hand, the lower value of $\mathcal{R}_0 (t)$ obtained
with the slower campaign ensures that in subsequent times this
strategy results in being more effective in lowering casualties.

\paragraph{Infinite Time Immunization:} Model~\eqref{eq:2} can
describe also the situation where the immunization provided by the
vaccine and/or acquired after recovering lasts for ever. In the case
$T_V \to +\infty$, system~\eqref{eq:2} becomes
\begin{equation}
  \label{eq:20}
  \left\{
    \begin{array}{l@{}}
      \displaystyle
      \dot S
      =
      - \rho_S \, I \, S
      + R (t, T_R)
      - p (t, S, V, I, R)
      \\
      \displaystyle
      \partial_t V + \partial_\tau V
      =
      - \rho_V \, V \, I
      \hfill \tau \in \mathopen[0, +\infty\mathclose[
      \\
      \displaystyle
      \dot I
      =
      \left(
      \rho_S \, S
      + \int_0^{+\infty} \!\! \rho_V (\tau) \, V (t,\tau) \d\tau
      + \int_0^{T_R}  \!\!\rho_R (\tau) \, R (t,\tau) \d\tau
      - (\theta + \mu)
      \right) I \qquad
      \\
      \displaystyle
      \partial_t R + \partial_\tau R
      =
      - \rho_R \, R \, I
      \hfill \tau \in [0, T_R]
      \\
      \displaystyle
      V (t, 0) = p(t,S, V, I, R)
      \\
      \displaystyle
      R (t, 0) = \theta \, I
    \end{array}
  \right.
\end{equation}
where it is clear that individuals that entered the $V$ population
will remain therein. An entirely similar system can be used to
describe the case $T_R \to +\infty$.

\begin{figure}[!h]
  \includegraphics[width = 0.5\linewidth]%
  {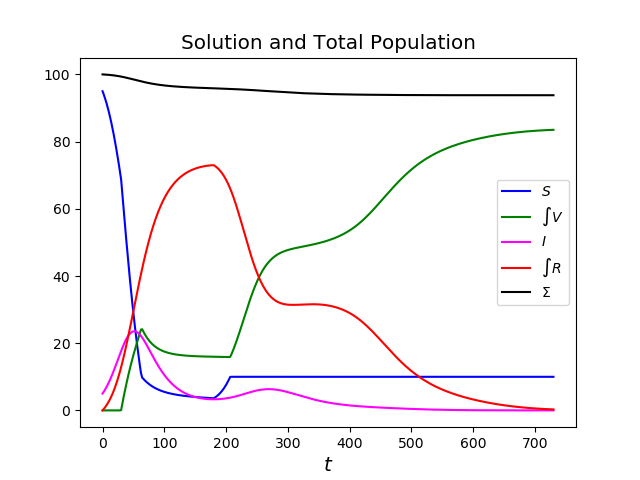}\hfil%
  \includegraphics[width = 0.5\linewidth]%
  {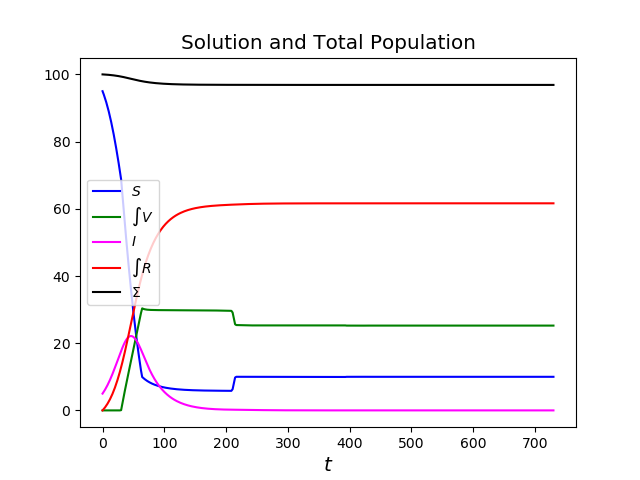}\\
  \caption{\small In these integration we used
    parameters~\eqref{eq:25}, the data~\eqref{eq:26} and the
    strategy~\eqref{eq:34} with $S_*=10$. Left, the case~\eqref{eq:20}
    where $T_V \to +\infty$ with $T_R = 180$.  Right, the case
    $T_R\to +\infty$ with $T_V = 180$. Note that epidemic waves
    essentially disappeared, quite quickly in the case
    $T_R \to +\infty$ on the right.\label{fig:INFI}}
\end{figure}

\begin{figure}[!h]
  \includegraphics[width = 0.5\linewidth]%
  {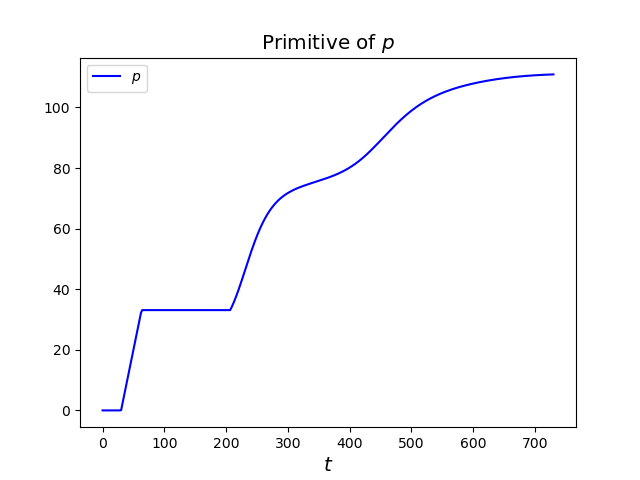}\hfil%
  \includegraphics[width = 0.5\linewidth]%
  {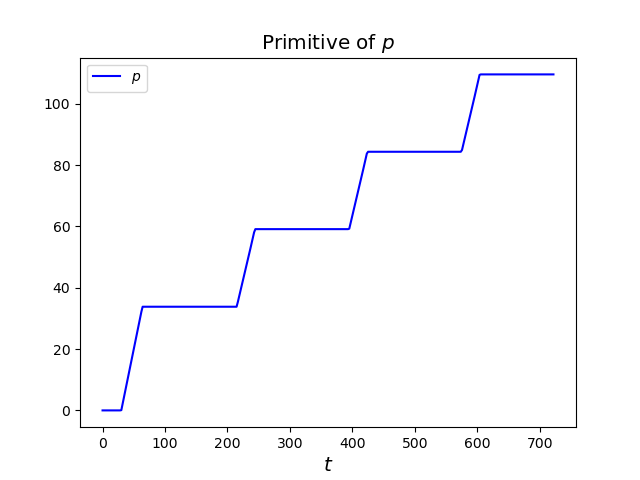}\\
  \caption{\small Distribution function
    $t \mapsto \int_0^t p (\tau) \d\tau$ of the vaccination strategies
    $p = p (t)$ in the integrations corresponding to
    parameters~\eqref{eq:25}, data~\eqref{eq:26} and
    strategy~\eqref{eq:34} with $S_*=10$. Left, the case~\eqref{eq:20}
    where $T_V \to +\infty$ with $T_R = 180$.  Right, the case
    $T_R\to +\infty$ with $T_V = 180$. Note that, on the right, the
    vaccination campaign is activated as soon as $S$ individuals are
    available, i.e., when vaccinated get back to be susceptible. 
    \label{fig:INFIp}}
\end{figure}

In the integrations below, we keep using the choice~\eqref{eq:25}, the
data~\eqref{eq:26} and the strategy~\eqref{eq:34} with $S_*=10$. The
resulting integrations, displayed in Figure~\ref{fig:INFI}, show an
evident stabilization effect induced by the infinite duration of the
immunization.
Indeed, epidemic waves are rather quickly smeared out, in particular
in the case $T_R \to +\infty$. As soon as individuals enter the $R$
population, they will (almost) never leave it, while all susceptible
individuals are vaccinated as soon as the effect of the previous
vaccination disappears, see Figure~\ref{fig:INFIp} on the right.

\section{Concurrent Vaccines}
\label{sec:Two}

We now consider the case of a vaccination campaign based on the
concurrent use of two different vaccines, so that vaccinated
individuals enter the population $V_1$ or $V_2$ depending on whether
they were dosed with vaccine $1$ or with vaccine $2$. Thus,
$V_1 (t,\tau)$, respectively $V_2 (t,\tau)$, measures the amount of
individuals at time $t$ that were dosed at time $t-\tau$ with vaccine
$1$, respectively with vaccine $2$. We also need to introduce the
controls specifying the speed at which the $2$ vaccines are dosed. The
equation for the $S$ population then reads:
\begin{equation}
  \label{eq:6}
  \begin{array}{rcl}
    \dot S
    & =
    & -\rho_S \, I \, S
      + V_1 (t,T_1)
      + V_2 (t,T_2)
      + R (t, T_R)
    \\
    &
    & - p_1 (t, S, V_1, V_2, I, R)
      - p_2 (t, S, V_1, V_2, I, R)
  \end{array}
\end{equation}
where we used obvious modification of the notation
in~\eqref{eq:2}. The above equation also prescribes that at time
$T_1$, individuals in the $V_1$ population get back to be susceptible,
and similarly for $T_2$.

Extending~\eqref{eq:2}, for the $V_1$, $V_2$ and $R$ populations we
obtain
\begin{equation}
  \label{eq:9}
  \left\{
    \begin{array}{@{}l@{\qquad}l@{}}
      \displaystyle
      \partial_t V_1 + \partial_\tau V_1
      =
      -\rho_1 (\tau) \, V_1 (t,\tau) \, I
      & \tau \in [0, T_1]
      \\
      \displaystyle
      \partial_t V_2 + \partial_\tau V_2
      =
      -\rho_2 (\tau) \, V_2 (t,\tau) \, I
      & \tau \in [0, T_2]
      \\
      \displaystyle
      \partial_tR + \partial_\tau R
      =
      -\rho_R (\tau) \, R (t,\tau) \, I
      & \tau \in [0, T_R]
      \\
      \displaystyle
      V_1 (t, 0) = p_1 (t, S, V_1, V_2, I, R)
      \\
      \displaystyle
      V_2 (t, 0) = p_2 (t, S, V_1, V_2, I, R)
      \\
      \displaystyle
      R (t, 0) = \theta \, I
    \end{array}
  \right.
\end{equation}
where, similarly to the previous section, the three time scales
$[0,T_1]$, $[0,T_2]$ and $[0, T_R]$ are entirely independent.

Finally, the $I$ population varies partly due to the propagation of
the infection and partly due to infected individuals recovering or
dying:
\begin{equation}
  \label{eq:10}
  \begin{array}{rcl}
    \dot I
    & =
    & \displaystyle
      \rho_S \, S \, I
      + \int_0^{T_1} \rho_1 (\tau) \, V_1 (t,\tau) \d\tau \, I
      + \int_0^{T_2} \rho_2 (\tau) \, V_2 (t,\tau) \d\tau \, I
    \\
    &
    & \displaystyle
      +  \int_0^{T_R} \rho_R (\tau) \, R (t,\tau) \d\tau \, I
      - \theta \, I - \mu \, I \,.
  \end{array}
\end{equation}

The natural extension of~\eqref{eq:6}--\eqref{eq:9}--\eqref{eq:10}
when $k$ different vaccines are available reads
\begin{equation}
  \label{eq:1}
  \left\{
    \begin{array}{l}
      \displaystyle
      \dot S = -\rho_S \, I \, S
      + \sum_{i=1}^k V_i (t,T_i)
      + R (t, T_R)
      - \sum_{i=1}^kp_i (t, S, V, I, R)
      \\
      \displaystyle
      \partial_tV_i+ \partial_\tau V_i
      =
      -\rho_i (\tau) \, V_i (t,\tau) \, I
      \hfill \tau \in [0, T_i]
      \qquad i=1 ,\ldots, k
      \\
      \displaystyle
      \dot I = \rho_S \, I \, S
      +
      \left(
      \sum_{i=1}^k \int_0^{T_i} \rho_i (\tau) \, V_i (t,\tau) \d\tau
      +
      \int_0^{T_R} \rho_R (\tau) \, R (t,\tau) \d\tau
      \right)
      I
      - \theta \, I
      - \mu \, I
      \\
      \displaystyle
      \partial_tR + \partial_\tau R
      =
      -\rho_R (\tau) \, R (t,\tau) \, I
      \hfill \tau \in [0, T_R]
      \\
      \displaystyle
      V_i (t, 0) = p_i (t, S, V, I, R)
      \\
      \displaystyle
      R (t, 0) = \theta \, I \,.
    \end{array}
  \right.
\end{equation}
Also in this general case, the index $\mathcal{R}_0 (t)$ can be
defined as
\begin{equation}
  \label{eq:1bis}
  \mathcal{R}_0 (t)
  =
  \dfrac{\rho_S \, S
    + \sum_{i=1}^k \int_0^{T_i} \!\! \rho_i (\tau) \, V_i (t,\tau) \d\tau
    + \int_0^{T_R}  \!\!\rho_R (\tau) \, R (t,\tau) \d\tau}{\theta+\mu}
\end{equation}
and identifies the times where $\dot I$ is positive or negative,
without explicitly requiring knowledge of $I (t)$.

\subsection{Comparing the Effects of Different Vaccines}
\label{subsec:comp-use-diff}

\begin{figure}[!h]
  \centering
  \begin{minipage}{0.5\linewidth}
    \includegraphics[width = \linewidth]{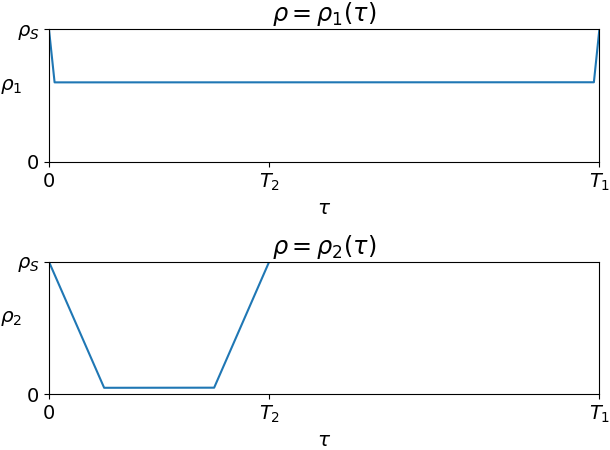}
  \end{minipage}%
  \begin{minipage}{0.5\linewidth}
    \begin{equation}
      \label{eq:7bis}
      \begin{array}{r@{\,}c@{\,}lr@{\,}c@{\,}l@{}}
        \rho_1 (\tau)
        & =
        & 6.0 \times 10^{-6}
        & \quad\mbox{for}\quad \tau
        & \in
        & [3, 297]
        \\
        T_1
        & =
        & 300
        \\[60pt]
        \rho_2 (\tau)
        & =
        & 5.0 \times 10^{-7}
        & \quad\mbox{for}\quad\tau
        & \in
        & [30, 90]
        \\
        T_2
        & =
        & 120
        \\[40pt]
      \end{array}
    \end{equation}
  \end{minipage}
  \caption{\small Characterizations of the vaccines used in the
    integrations of~\eqref{eq:1} with $k=2$ and the
    choices~\eqref{eq:7}--\eqref{eq:2}--\eqref{eq:25}--\eqref{eq:26}.
    Vaccine~$1$ provides a weak protection almost immediately, while
    vaccine~$2$ is more protective but after a while and for a shorter
    time. \label{fig:cinque}}
\end{figure}

In the integrations of this paragraph we keep using the
choices~\eqref{eq:2}--\eqref{eq:25}--\eqref{eq:26}, so that the
reference situation with no vaccination campaign is the one discussed
in~\S~\ref{subsec:comp1} and illustrated in
Figure~\ref{fig:no0Vacc}. We introduce $2$ vaccines, say $1$ and $2$,
characterized by the diagrams in Figure~\ref{fig:cinque}, see
also~\eqref{eq:7bis}.

We consider the strategies
\begin{enumerate}[label=\textbf{(1\arabic*):},
  ref=\textbf{\textup{(1\arabic*)}}, align=left]

\item \label{case:11} Vaccine~$1$ is used throughout, from time $t=30$
  on, while Vaccine~$2$ is not used.

\item \label{case:12} Vaccine~$1$ is used for $t \in [30, 380]$, while
  Vaccine~$2$ is used for $t \in [380, 730]$.
\end{enumerate}
\vspace{-.4cm}
\begin{enumerate}[label=\textbf{(2\arabic*):},
  ref=\textbf{\textup{(2\arabic*)}}, align=left]

\item \label{case:21} Vaccine~$2$ is used for $t \in [30, 380]$, while
  Vaccine~$1$ is used for $t \in [380, 730]$.

\item \label{case:22} Vaccine~$2$ is used throughout, from time $t=30$
  on, while Vaccine~$1$ is not used.
\end{enumerate}
\vspace{-.5cm}
\begin{enumerate}[label=\textbf{(1/2):}, ref=\textbf{\textup{(1/2)}},
  align=left]

\item \label{case:1-over-2} Both vaccines are used throughout, with
  the same number of doses.
\end{enumerate}

\noindent In the present setting we are assuming that vaccinated
individuals that get back to be susceptibles are vaccinated as soon as
possible. Therefore, it is intuitive that Vaccine~$2$ results in being
the best choice, as shown in Table~\ref{tab:cinque}.
\begin{table}[!h]
  \centering
  \begin{tabular}[c]{l|r|r|r|r|r|r}
    Strategy
    & ref.
    & (11)
    & (12)
    & (21)
    & (22)
    & (1/2)
    \\ \hline
    Deaths
    & 15.1
    & 11.4
    & 7.20
    & 4.26
    & 3.92
    & 5.78
    \\ \hline
    1 Doses
    & 0.00
    & 187
    & 125
    & 131
    & 0.00
    & 159
    \\
    2 Doses
    & 0.00
    & 0.00
    & 151
    & 215
    & 480
    & 159
    \\
    Doses Tot.
    & 0.00
    & 187
    & 275
    & 346
    & 480
    & 318
  \end{tabular}
  \caption{\small Statistics on the solutions to~\eqref{eq:1} with $k=2$ and
    the choices~\eqref{eq:25}--\eqref{eq:26}--\eqref{eq:7bis} corresponding to the strategies outlined in~\S~\ref{subsec:comp-use-diff}. The leftmost column refers to the reference solution where no vaccination takes place, see~\S~\ref{subsec:comp1}.
    \label{tab:cinque}}
\end{table}
Indeed, in the present framework, once an individual is vaccinated
with a Vaccine~$1$, he/she can not be vaccinated using the more
efficient Vaccine~$2$ as long as the first immunization is, though
only poorly, effective. This also explains the different outcomes of
the strategies~\ref{case:12} and~\ref{case:21}.  Note also that
strategy~\ref{case:22} allows to dose $80\%$ of the initial population
$5$ times, see Table~\ref{tab:cinque}.

\section{Continuous and Discrete Age Structures}
\label{sec:Gen}

Age differences can play a significant role in the reaction of
individuals to the infection. We thus extend our framework to account
also for age differences. First, we insert a continuous age structure,
later in \S~\ref{subsec:age-depend-vacc} we consider discrete age
classes. In the first case $a \in \reali_+$ is a continuous variable
and the convective terms
$\partial_a S, \partial_a V, \partial_a I, \partial_a R$
in~\eqref{eq:11} describe the aging of the individuals in the
population $S,V,I,R$. In the latter case, $a$ is a discrete variable
ranging in the finite set of the age classes considered and no aging
term is present, see~\eqref{eq:14}. The former approach seems more
accurate, but on short time intervals the second is a usual and
acceptable simplification.

For simplicity, we detail the age structured version of~\eqref{eq:2}
corresponding to only one vaccine. The extension of the $k$ vaccines
case~\eqref{eq:1} being only technically more intricate. We thus
obtain:
\begin{equation}
  \label{eq:11}
  \left\{
    \begin{array}{l@{}}
      \displaystyle
      \partial_t S + \partial_a S
      =
      - \rho_S \, I \, S
      + V \left(t, a, T_V (a)\right)
      + R \left(t, a, T_R (a)\right)
      - p (t, a, S, V, I, R)
      \\
      \displaystyle
      \partial_t V  + \partial_a V + \partial_\tau V
      =
      - \rho_V \, V \, I
      \hfill \tau \in [0, T_V (a)]
      \\
      \displaystyle
      \partial_t I  + \partial_a I
      =
      \rho_S \, S \, I
      + \int_0^{T_V} \!\! \rho_V (a, \tau) \, V (t, a, \tau) \d\tau \, I
      + \int_0^{T_R}  \!\!\rho_R (a, \tau) \, R (t, a, \tau) \d\tau \, I
      \\
      \qquad\qquad\qquad- \theta \, I - \mu \, I
      \\
      \displaystyle
      \partial_t R  + \partial_a  R + \partial_\tau R
      =
      - \rho_R \, R \, I
      \hfill \tau \in [0, T_R (a)]
      \\
      \displaystyle
      V (t, a, 0) = p(t,a, S, V, I, R)
      \\
      \displaystyle
      R (t, a, 0) = \theta \, I \,.
    \end{array}
  \right.
\end{equation}
Note that here all effects of vaccines are age dependent. The
immunization time provided by the vaccine is $T_V = T_V (a)$ and,
similarly, also the immunization ensured by recovering from the
disease is age dependent: $T_R = T_R (a)$. Remark that~\eqref{eq:11}
is able to take into consideration the different effectiveness of the
vaccine at different ages, thanks to the dependence of $\rho_V$ also
on $a$: $\rho_V = \rho_V (a, \tau)$. Similarly, in~\eqref{eq:11} also
the recovery rate $\theta$ depends on the age, i.e.,
$\theta = \theta(a)$, as well as the mortality rate, $\mu = \mu (a)$.

As usual in age structured models, further boundary conditions need to
be supplemented, taking care of the newborns, such as
\begin{equation}
  \label{eq:12}
  S (t, 0) = b (t)
  \,,\qquad
  V (t,0,\tau) = 0
  \,,\qquad
  I (t,0) = 0
  \,,\qquad
  R (t,0,\tau) = 0 \,,
\end{equation}
where $b (t)$ is the time dependent natality. Other natality terms can
be considered, depending, for instance, on the total amount of
susceptibles.

However, typically, the use of a pandemic model may be of interest on
time intervals far smaller than the average life span of
individuals. Therefore, it is convenient to consider a fixed number,
say $m$, of age classes. As a consequence, we have $m$ different
populations of susceptibles, of vaccinated, infected and recovered,
obtaining the mixed multiscale system
\begin{equation}
  \label{eq:14}
  \left\{
    \begin{array}{l}
      \displaystyle
      \dot S_a
      =
      -  \left(\sum_{\alpha=1}^m \rho_S^{a,\alpha} \, I_\alpha\right) S_a
      + V_a (t,T_V^a)
      + R_a (t, T_R^a)
      - p_a (t, S, V, I, R)
      \\
      \displaystyle
      \partial_t V_a + \partial_\tau V_a
      =
      - \left(\sum_{\alpha=1}^m \rho_V^{a,\alpha} \, I_\alpha\right) V_a
      \hfill \tau \in [0, T_V^a]
      \\
      \displaystyle
      \dot I_a
      =
      \left(\sum_{\alpha=1}^m \rho_S^{a,\alpha} \, I_\alpha\right) S_a
      +
      \left(\sum_{\alpha=1}^m I_\alpha \int_0^{T_V^a} \rho_V^{a,\alpha} (\tau) \, V_a (t,\tau) \d\tau \right)
      \\
      \qquad\qquad
      \displaystyle
      +
      \left(\sum_{\alpha=1}^m I_\alpha \int_0^{T_R^a} \rho_R^{a,\alpha} (\tau) \, R_a (t,\tau) \d\tau \right)
      - \theta_a \, I_a
      - \mu_a \, I_a
      \\
      \displaystyle
      \partial_t R_a + \partial_\tau R_a
      =
      - \left(\sum_{\alpha=1}^m \rho_R^{a,\alpha} \, I_\alpha\right) R_a
      \hfill \tau \in [0, T_R^a]
      \\
      \displaystyle
      V_a (t, 0) = p_a(t,S, V, I, R)
      \\
      \displaystyle
      R_a (t, 0) = \theta_a \, I_a
    \end{array}
  \right.
  a=1, \ldots, m\,.
\end{equation}
Above, the terms $\rho_S^{a,\alpha}$, $\rho_V^{a,\alpha}$ and
$\rho_R^{a,\alpha}$ quantify the spreading of the virus between the
age class $a$ and the age class $\alpha$ in the populations $S$, $V$
and $R$.

In~\eqref{eq:14}, differently from what happens in~\eqref{eq:11}, the
total number of individuals in the age class $a$ may vary only due to
the mortality in that class, i.e.,
\begin{equation}
  \label{eq:13}
  \dfrac{\d{~}}{\d{t}}
  \left(
    S_a (t)
    +
    \int_0^{T_V^a} V_a (t,\tau)  \d\tau{}
    +
    I_a (t)
    +
    \int_0^{T_R^a} R_a (t,\tau)  \d\tau{}
  \right)
  =
  -\mu_a \, I_a (t) \,,
\end{equation}
so that infection propagates among individuals of different classes,
but no individual changes its age class.

The introduction of an index similar to $\mathcal{R}_0 (t)$ is
formally possible, but the resulting expression necessarily explicitly
depends on $I_a (t)$.

\subsection{Comparing Age Dependent Vaccination Strategies}
\label{subsec:age-depend-vacc}

We limit the numerical integrations of~\eqref{eq:14} to the case of
only $2$ classes, say the \emph{young} one (indexed with $1$) and the
\emph{old} one ($2$). For a different approach to the modeling of $2$
age classes, refer for instance to~\cite{VerrelliDellaRosa}.

\paragraph{The Reference Situation:} Consider first the case where no
vaccination campaign takes place. On the basis of a qualitative
approach as in \S~\ref{subsec:comp1}, we choose the following set of
parameters:
\begin{equation}
  \label{eq:15}
  \begin{array}{r@{\,}c@{\,}l@{\qquad}%
    r@{\,}c@{\,}l@{\qquad}r@{\,}c@{\,}l@{\qquad}r@{\,}c@{\,}l}
    \rho_S^{11}
    & =
    & 3.0 \times 10^{-3}
    & \rho_S^{12}
    & =
    & 1.0\times 10^{-3}
    & \rho_S^{21}
    & =
    & 2.0 \times 10^{-3}
    & \rho_S^{22}
    & =
    & 1.0\times 10^{-3}
    \\
    \theta^1
    & =
    & 6.0\times10^{-2}
    & \theta^2
    & =
    & 4.0\times10^{-2}
    & \mu^1
    & =
    & 5.0 \times 10^{-4}
    & \mu^2
    & =
    & 2.0\times10^{-3}
    \\
    T_R^1
    & =
    & 180
    & T_R^2
    & =
    & 140
    &\rho_R^-
    & =
    & 2.0\times 10^{-5}
  \end{array}
\end{equation}
and we keep referring to the choices of $\Phi$ and $\Psi$
in~\eqref{eq:7} and~\eqref{eq:8}, so that for $a,\alpha = 1,2$
\begin{equation}
  \label{eq:16}
  \rho_V^{a,\alpha} (\tau)
  =
  \rho_V^- + (\rho_S^{a,\alpha}-\rho_V^-) \;
  \Psi (\tau / T_V^{a,\alpha})
  \quad \mbox{ and } \quad
  \rho_R^{a,\alpha} (\tau)
  =
  \rho_R^- + (\rho_S^{a,\alpha}-\rho_R^-) \;
  \Phi (\tau / T_R^{a,\alpha}) \,.
\end{equation}
The above choices reflect the fact that class $2$ individuals suffer a
higher mortality ($\mu_2 = 4\, \mu_1$) and have a slower recovery
($\theta^2= 0.67\, \theta^1$). The two age classes differ also in the
time scales, the younger ones having longer periods of (partial)
immunization both after recovery and after vaccination. On the other
hand, among class $1$ individuals the virus spreads faster
($\rho^{11}_S / \rho^{22}_S = 3$).

Throughout, we carry the integrations up to a final time $730$
(roughly corresponding to $2$ years) and with the initial datum (for
$a=1,2$)
\begin{equation}
  \label{eq:17}
  S_1 (0) = 42\,,\quad
  S_2 (0) = 53\,,\quad
  V_a (0,\tau) = 0 \,,\quad
  I_1 (0) = 1\,,\quad
  I_2 (0) = 4\,,\quad
  R_a (0,\tau) = 0 \,.
\end{equation}

\begin{figure}[!h]
  \includegraphics[width = 0.33\linewidth, trim=20 30 20 20]%
  {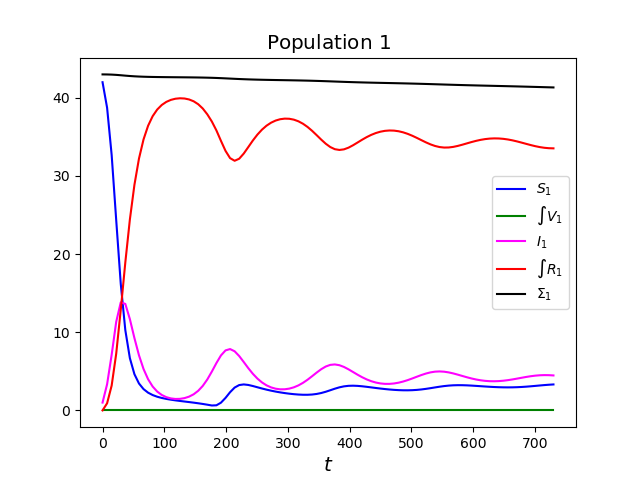}\hfil%
  \includegraphics[width = 0.33\linewidth, trim=20 30 20 20]%
  {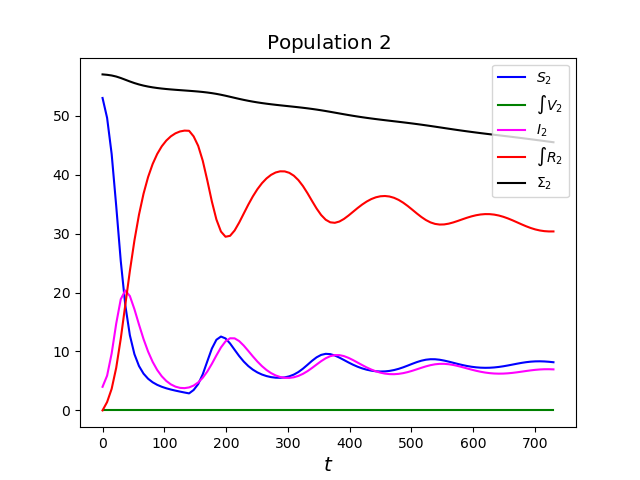}\hfil%
  \includegraphics[width = 0.33\linewidth, trim=20 30 20 20]%
  {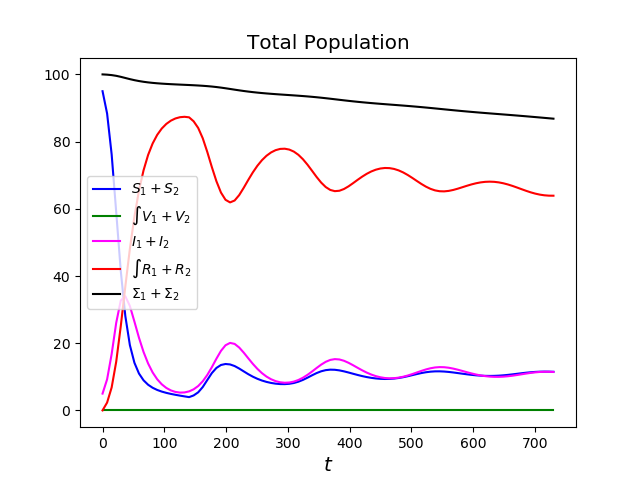}\\
  \caption{\small Solution to~\eqref{eq:14} with
    parameters~\eqref{eq:15}--\eqref{eq:16} and initial
    datum~\eqref{eq:17} in the case no vaccination is dosed. Left,
    population~$1$; middle, population~$2$ and, right, their sum. In
    all diagrams, $\Sigma$ stands for the sum of all individuals of
    the considered age class(es).\label{fig:vageq0}}
\end{figure}

The resulting evolution is displayed in Figure~\ref{fig:vageq0}.
In this reference situation, the casualties after time $730$ (i.e.,
$2$ years) amount to $1.67$ in class~$1$, $11.5$ in class~$2$,
totaling to~$13.2$ (the total initial population being $100$). Note
the formation of \emph{``epidemic waves''}~\cite{LemonMahmoud2005},
Figure~\ref{fig:vageq0}.

In this framework, a variety of age dependent vaccination strategies
can be adopted. Concerning the vaccine, following~\eqref{eq:25}, we
keep the following choices fixed:
\begin{equation}
  \label{eq:19}
  \rho_V^{a,\alpha} (\tau) = \rho_V^- + (\rho_S^{a,\alpha}-\rho_V^-) \; \Psi (\tau/T_V^a)
  \mbox{ for } a,\alpha=1,2 \qquad \mbox{ and } \quad
  \begin{array}{r@{\;}c@{\;}l}
    T_V^1
    & =
    & 200
    \\
    T_V^2
    & =
    & 160
    \\
    \rho_V^-
    & =
    & 1.0 \times 10^{-5}
  \end{array}
\end{equation}
where $\Psi$ is as in~\eqref{eq:8}. Throughout, we let the vaccination
campaign begin after time $t=30$. We consider below $4$ instances:
\begin{enumerate}[label=\textbf{Feedback:},
  ref=\textbf{\textup{Feedback}}, align=left]

\item \label{case:feedback} $p_a (t)$ is proportional to the number of
  infected individuals in class~$a$, i.e.,
  $p_a (t) = I_a (t) / \left(I_1 (t) + I_2 (t)\right)$ as long as
  $S_a (t) > 0$, for $a=1,2$.
\end{enumerate}
\begin{enumerate}[label=\textbf{Half--Half:},
  ref=\textbf{\textup{Half--Half}}, align=left]

\item \label{case:half-half} $p_a (t) = 0.5$ as long as there are
  susceptibles in class $a$, i.e., $S_a (t)>0$, for $a=1,2$.
\end{enumerate}
\begin{enumerate}[label=\textbf{Class~$2$ First:},
  ref=\textbf{\textup{Class~$2$ First}}, align=left]

\item \label{case:2-first} for $t \in [30, 380]$, $p_1 (t) = 0$ and
  $p_2 (t) = 1$ as long as $S_2 (t)>0$; for $t \in [380, 730]$,
  $p_1 (t) = 1$ and $p_2 (t) = 0$ as long as $S_1 (t)>0$.
\end{enumerate}
\begin{enumerate}[label=\textbf{Class~$1$ First:},
  ref=\textbf{\textup{Class~$1$ First}}, align=left]

\item \label{case:1-first} for $t \in [30, 380]$, $p_1 (t) = 1$ and
  $p_2 (t) = 0$ as long as $S_1 (t)>0$; for $t \in [380, 730]$,
  $p_1 (t) = 0$ and $p_2 (t) = 1$ as long as $S_2 (t)>0$.
\end{enumerate}

\noindent In all strategies, the total number of vaccines dosed per
day is at most $1\%$ of the total initial population, as soon as the
number of susceptibles is sufficiently high. In these examples, we
also let vaccinated be dosed again as soon as they get back to be
susceptible. The present framework clearly allows also to leave an
amount of non vaccinated individuals, as in~\S~\ref{subsec:comp1}.

It is evident that \ref{case:1-first} is likely to be the least
effective strategy, as it actually results from
Table~\ref{tab:q}. Less intuitive is the fact that \ref{case:2-first}
is only slightly better, in particular comparing the total number of
casualties.
\begin{table}[!h]
  \centering
  \begin{tabular}[c]{l|r|r|r|r|r}
    Strategy
    & Reference
    & Feedback
    & Half--Half
    & Class $2$ First
    & Class $1$ First
    \\ \hline
    $1$ Deaths
    & 1.67
    & 0.620
    & 0.630
    & 1.24
    & 1.35
    \\
    $2$ Deaths
    & 11.5
    & 3.52
    & 3.66
    & 7.47
    & 8.51
    \\
    Deaths Tot.
    & 13.2
    & 4.14
    & 4.29
    & 8.71
    & 9.86
    \\ \hline
    $1$ Doses
    & 0.00
    & 106
    & 106
    & 41.8
    & 28.5
    \\
    $2$ Doses
    & 0.00
    & 206
    & 203
    & 87.0
    & 84.4
    \\
    Doses Tot.
    & 0.00
    & 313
    & 309
    & 129
    & 113
  \end{tabular}
  \caption{\small Statistics on the solutions to~\eqref{eq:14} with
    parameters~\eqref{eq:15}--\eqref{eq:16}, initial
    datum~\eqref{eq:17} and with the vaccination strategies detailed
    in~\S~\ref{subsec:age-depend-vacc}.\label{tab:q}}
\end{table}
Surprisingly, \ref{case:2-first} results in a number of casualties in
class~$1$ even lower than that resulting from strategy
\ref{case:1-first}.
\begin{figure}[!h]
  \includegraphics[width=0.5\linewidth]{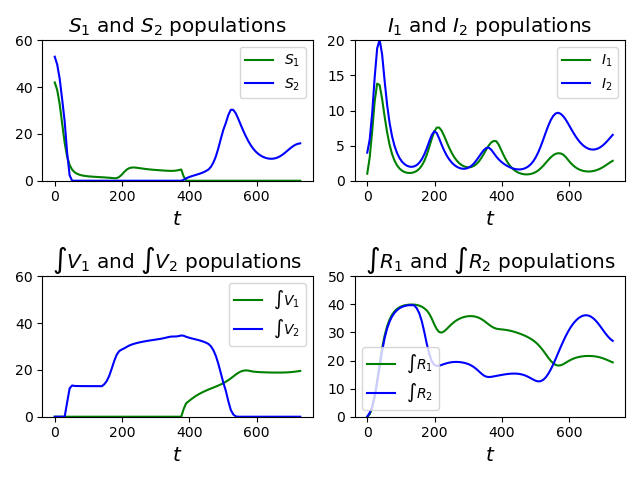}%
  \includegraphics[width=0.5\linewidth]{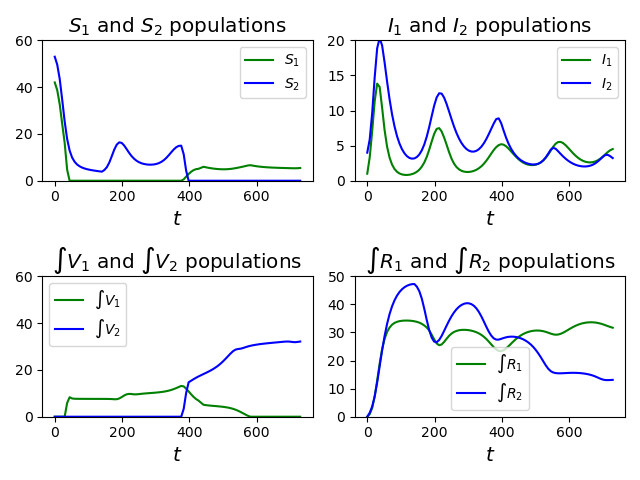}%
  \caption{\small Graphs of the solutions to~\eqref{eq:14} with
    parameters~\eqref{eq:15}--\eqref{eq:16}, initial
    datum~\eqref{eq:17} and with the vaccination strategies
    \ref{case:2-first}, left, and \ref{case:1-first}, right, as
    detailed in~\S~\ref{subsec:age-depend-vacc}.\label{fig:q2}}
\end{figure}
Moreover, in both \ref{case:1-first} and \ref{case:2-first}
strategies, the rise of rather persisting epidemic waves is
evident. In particular, in the latter case the final increase in the
number of infected of both classes induces to expect a worsening of
the situation in the long run.

Among the strategies considered, the one resulting most effective in
containing casualties is the \ref{case:feedback} one. However, it is
not easy to anticipate that, mainly due to the particular initial data
chosen, it is only slightly better than the \ref{case:half-half}
one. Indeed, a feedback strategy is generally prone to provide better
results than an open loop one.

\begin{figure}[!h]
  \includegraphics[width=0.5\linewidth]{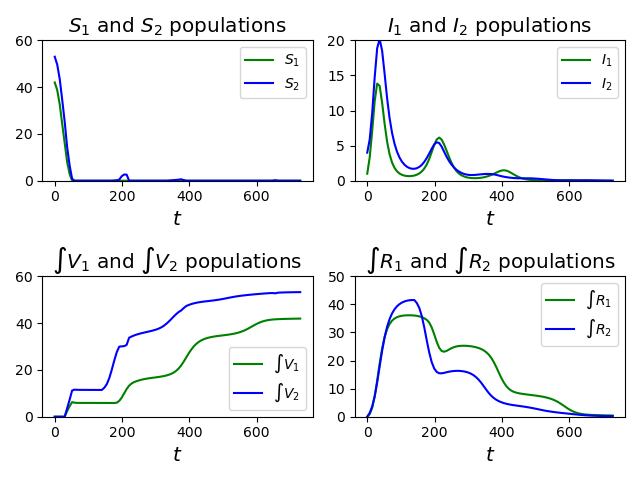}%
  \includegraphics[width=0.5\linewidth]{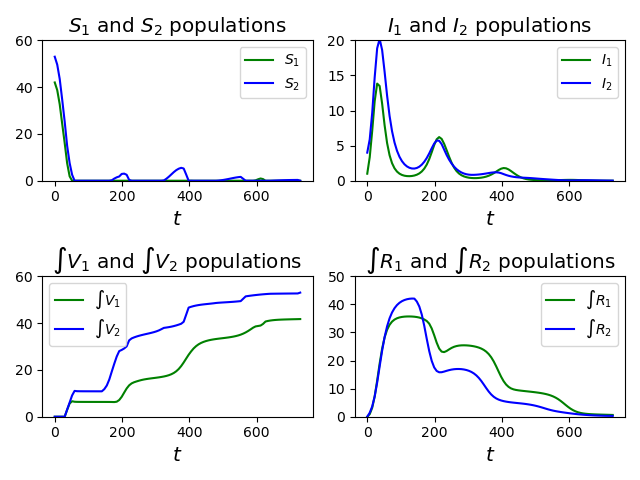}%
  \caption{\small Graphs of the solutions to~\eqref{eq:14} with
    parameters~\eqref{eq:15}--\eqref{eq:16}, initial
    datum~\eqref{eq:17} and with the vaccination strategies
    \ref{case:feedback}, left, and \ref{case:half-half}, right, as
    detailed in~\S~\ref{subsec:age-depend-vacc}.\label{fig:q1}}
\end{figure}

It stems out of these examples that, in order to reduce the number of
casualties, it is of key importance to bound the number of
susceptibles, as a comparison between the graphs in
Figure~\ref{fig:q2} and in Figure~\ref{fig:q1} shows.

It is also worth noting that a \emph{``weak''} vaccination campaign
can lead to somewhat persistent epidemic waves. Indeed, compare the
qualitative behavior of the maps $t \mapsto I_a (t)$ in the reference
case (Figure~\ref{fig:vageq0}), in the successful cases
\ref{case:feedback} or \ref{case:half-half} to those corresponding to
the strategies \ref{case:1-first} or \ref{case:2-first}
(Figure~\ref{fig:q1}). The epidemic waves in the latter case appear to
be quite persistent, while they fade out sooner in the former 2
cases. Indeed, the lack of any vaccination results in a high mortality
that hinders the repeated formation of waves.
On the opposite, an efficient vaccination campaign quickly flattens
the $S_a$ curve near to $0$. A weak vaccination campaign still reduces
the number of casualties but may be not sufficiently strong to
eradicate the disease, which thus keeps propagating in waves.

\section{Parameters' Choices in the Case of Covid--19}
\label{sec:real-data-fitting}

Here we present and justify possible \emph{a priori} choices for the
parameters and functions entering our framework on the basis of
available measurements. A different approach might consist in fitting,
\emph{a posteriori}, the solutions to the above models to the measured
evolution.

As is well known, the actual values of the parameters or functions
$\rho_S, \rho_V, \rho_R, \theta, \mu$ appearing in the various models
depend on possible normalizations of the total number of
individuals. In the previous sections, for instance, we set the total
initial population to $100$. The time scales, defined for instance by
$T_V$ and $T_R$ in the case of~\eqref{eq:2}, can be directly deduced
from the literature.

For simplicity, we refer to~\eqref{eq:2}, the further parameters
entering the other models can be evaluated similarly.

\paragraph{Parameter $\rho_S$: } In any attempt to obtain real values
from a predictive model, this parameter has to be taken as time
dependent. Indeed, not only it heavily depends on the introduction of
lockdown restrictions or on the season, the mere social awareness of
the disease may significantly affect its value. As an example, we
refer to~\cite[Table~3]{Law2020}, where time dependent values for
$\rho_S$ (there denoted $z\, \beta_t$) are deduced from real data.

\paragraph{Parameter $T_V$ and Function $\rho_V = \rho_V (\tau)$: }
Here, by vaccination we intend the full treatment consisting of $2$
injections, dosed sufficiently near so that there is no loss in the
protection they provide. Therefore, we rely on a function $\rho_V$ of
the type in Figure~\ref{fig:rhoV}, with $T_V$ being the duration of
the immunization provided by the $2$ doses. Due to the relatively
short history Covid--19 vaccinations, this datum can only be inferred,
see for instance~\cite{Baden2021403, Polack20202603}. On the basis
of~\cite{Zenilman20211761, Whitley20211774}, it seems safe to assume
$T_V = 180 \mbox{ days}$, i.e., $6$ months, for the BNT162B2 mRNA
vaccine as well as for the mRNA-1273 vaccine. A choice of
$\rho_V^- = 0.05 \, \rho_S$ is realistic in the case of the BNT162b2
mRNA vaccine, see \cite{Polack20202603}, while
$\rho_V^- = 0.06 \, \rho_S$ seems justified by~\cite{Baden2021403} in
the case of the mRNA-1273 vaccine.

\paragraph{Parameter $T_R$ and Function $\rho_R = \rho_R (\tau)$: }
The current literature provides different results about the duration
of the immunization enjoyed by those who recovered from Covid--19. For
instance, in~\cite{Ripperger2020925}, the authors suggest for $T_R $
a value between $150$ and $210$ days. Assuming for the function
$\rho_R = \rho_R (\tau)$ a shape like that in~\eqref{eq:7}, we are
left to estimate $\rho_R^-$, which measures the best level of
protection provided by the recovery. In the data collected
in~\cite{Shrestha2021.06.01.21258176}, none of the $1359$ non
vaccinated that recovered got infected. In~\cite{Lumley2021533}, out
of $11364$ individuals that recovered, only $2$ resulted infected
according to a PCR test. Clearly, it is natural to expect that both
$T_R$ and $\rho_R^-$ are significantly age dependent.

\paragraph{Parameters $\theta$ and $\mu$: } Both these parameters should
better be considered time dependent whenever simulations are meant to
provide results on a scale of several months. Indeed, care protocols
have been continuously updated since Covid--19 outbreak and new drugs
have been introduced. As a reference, we recall that in it is
suggested in~\cite{Russo20212517} that $\theta/\mu \approx 18.2$, obtained from
statistics in the Milan (Italy) area between February 19th, 2020 and
January 21st, 2021.

\paragraph{Function $p = p (t)$: } This function quantifies how many
vaccination are dosed per unit time (e.g.~day). Clearly, it is time
dependent and its value has been chosen according to different
policies in different nations. Often, health care workers were given
the highest priority with old or fragile individuals coming next. As
reference values, we record that in Italy on August 4th 2021, $171565$
individuals (i.e., about $0.29\%$ of the Italian population) received
their first dose, while they were $5544$ ($0.0093\%$) on November 1st,
2021, data taken from~\cite{covid19GitHub}.

\section{Conclusion}
\label{sec:conclusion}

This paper introduces a framework for the multiscale modeling of a
vaccination campaign in presence of a pandemic.  Different concurrent
vaccines can be considered, each of them is characterized by its own
efficiency and provides an immunization whose level depends on the
time since dosing.  Different age classes can be considered to account
for the dependence of mortality, infectivity, vaccine efficiency... on
age.  Within this framework, different vaccination strategies can be
simulated and compared, in terms of casualties, number of infected
individuals or number of vaccines dosed, for instance.

To our knowledge, a general theorem ensuring the well posedness of
these models is still unavailable.  We expect that this result is at
reach along the lines in~\cite{MauroSIR, futuro}.  Also the search for
an \emph{``optimal''} vaccination strategy is still an open problem,
however recent results are pointing in this direction; see for
instance to~\cite{DiGiamberardinoEtAl, keimerPflug, McQuadeEtAl}.

Aiming at quantitatively reliable forecasts by means of the present
framework requires accurate knowledge of various data. In particular,
the efficiency of vaccines, here quantified through the function
$\rho_V = \rho_V (\tau)$ (or $\rho_V = \rho_V (\tau,a)$), appears as quite difficult. In
this connection, it looks promising to deal with the uncertainties
intrinsic to these functions through the recent techniques
in~\cite{AlbiBertagliaEtAl, MR4263205, MR4313914}.

The present framework is quite flexible and several extensions are
easily at reach. For instance, letting $\rho_S$ and/or $\rho_R$ depend
explicitly also on time $t$ may account for the insurgence of new
virus mutations or strict lockdown policies.  Spatial movements can be
incorporated using exactly the same techniques as
in~\cite[Section~6]{FrancescaElena}
or~\cite[Formula~(1)]{preprint}. Gender differences only amount to
introduce further distinctions among the unknown variables.

\appendix
\section{Appendix: Proof of~\eqref{eq:29} and~\eqref{eq:30} }

\begin{lemma}
  \label{lem:uuu}
  The solution to $\:\left\{
    \begin{array}{@{\,}l}
      \partial_t u + \partial_\tau u = u \; \phi (t)
      \\
      u (t,0) = \alpha \, \phi (t)
    \end{array}
  \right.$ is $u (t,\tau) = \alpha \; \phi (t-\tau) \, e^{\int_{t-\tau}^t \phi (s)\d{s}}$.
\end{lemma}

\noindent The proof is a straightforward computation, hence it is
omitted.

\begin{lemma}
  \label{lem:si}
  Assume that, for a suitable $p$, problem~\eqref{eq:2} admits a
  classical solution on $\reali_+$ with $\mathcal{R}_0 (t)=r_*$ for
  all $t \in \left[t_*, +\infty\right[$, for a suitable $t_* > 0$ and
  with $\mathcal{R}_0$ as defined in~\eqref{eq:3}. Then, \eqref{eq:29}
  and~\eqref{eq:30} hold.
\end{lemma}

\begin{proof}
  The above assumptions ensure that $S,V,I,R$ also solve for
  $t \geq t_*$ the mixed ODE--PDE problem
  \begin{equation}
    \label{eq:32}
    \left\{
      \begin{array}{ll@{}}
        \displaystyle
        \dot S
        =
        - \rho_S \, I \, S
        + V (t,T_V)
        + R (t, T_R)
        - p (t, S, V, I, R)
        \\
        \displaystyle
        \partial_t V + \partial_\tau V
        =
        - \rho_V \, V \, I
        & \tau \in [0, T_V]
        \\
        \displaystyle
        \dot I
        =
        (r_*-1) \, (\theta+\mu) \,I
        \\
        \displaystyle
        \partial_t R + \partial_\tau R
        =
        - \rho_R \, R \, I
        & \tau \in [0, T_R]
        \\
        \displaystyle
        V (t, 0) = p(t,S, V, I, R)
        \\
        \displaystyle
        R (t, 0) = \theta \, I
      \end{array}
    \right.
  \end{equation}
  We then obtain a closed equation for $I$, namely
  $\dot I = (r_*-1) \, (\theta+\mu) \,I$, whose solution is the first
  line in~\eqref{eq:29}. Then, the initial--boundary value problem
  \begin{displaymath}
    \left\{
      \begin{array}{l}
        \displaystyle
        \partial_t R + \partial_\tau R
        =
        - \rho_R \, R \, I
        \\
        \displaystyle
        R (t, 0) = \theta \, I
      \end{array}
    \right.
  \end{displaymath}
  fits into Lemma~\ref{lem:uuu} with $u (t,\tau) = R (t,\tau)$,
  $\phi (t) = -\rho_R \, I (t)$, $\alpha = -\theta/\rho_R$, proving
  the second line in~\eqref{eq:29}.

  Verifying~\eqref{eq:30} is now straightforward.
\end{proof}

\section{Appendix: A Note on the Numerical Algorithm Adopted}
\label{sec:append-note-numer}

The systems considered consist of mixed Ordinary--Partial Differential
Equations. All differential equations are first order, non linear and
leave $\mathopen[0, +\infty\mathclose[$ invariant. The particular
structure of the convective parts in the PDEs, where both independent
variables are times, suggests to use a simple upwind
scheme~\cite[\S~4.2]{LeVequeBook2002}, using the same mesh $\Delta t$
for all independent variables ($t$ and $\tau$), although they vary in
different time interval. The right hand sides of all equations are
computed through a first order forward Euler method, taking care that
equality~\eqref{eq:27} (or its analog~\eqref{eq:13}) keeps holding at
each time step.

To prevent the $S$ variable getting negative when it is near to $0$,
we employed a simple predictor-corrector method. Whenever
$S (t+\Delta t)$ gets negative, the value of $p (t)$ is recomputed,
consistently with the equation, so that $S (t+\Delta t) =0$.

\section*{Acknowledgments}
The authors was partly supported by the GNAMPA~2020 project
\emph{"From Wellposedness to Game Theory in Conservation Laws"}. The
\emph{IBM Power Systems Academic Initiative} substantially contributed
to all numerical integrations.

\paragraph{Funding:} Not applicable.

\paragraph{Data Availability: } Not applicable.

\paragraph{Code Availability:} Not applicable.

\paragraph{Declarations:}

\subparagraph{Conflict of interest:} The authors declare that they have no conflict of interest.

{\small

  \bibliographystyle{abbrv}

  \bibliography{Ultima}

\begin{thebibliography}{10}

\bibitem{covid19GitHub}
{C}{O}{V}{I}{D}-19 opendata vaccini.
\newblock \url{https://github.com/italia/covid19-opendata-vaccini}.
\newblock Accessed: 2021-12-02.

\bibitem{Al-QanessEtAl}
M.~Al-Qaness, A.~Ewees, H.~Fan, and M.~Aziz.
\newblock Optimization method for forecasting confirmed cases of
  {C}{O}{V}{I}{D}-19 in {C}hina.
\newblock {\em Applied Sciences}, 9(3), 2020.

\bibitem{AlbiBertagliaEtAl}
G.~Albi, G.~Bertaglia, W.~Boscheri, G.~Dimarco, et~al.
\newblock Kinetic modelling of epidemic dynamics: social contacts, control with
  uncertain data, and multiscale spatial dynamics, 2021.

\bibitem{MR4263205}
G.~Albi, L.~Pareschi, and M.~Zanella.
\newblock Control with uncertain data of socially structured compartmental
  epidemic models.
\newblock {\em J. Math. Biol.}, 82(7):Paper No. 63, 41, 2021.

\bibitem{MR4313914}
G.~Albi, L.~Pareschi, and M.~Zanella.
\newblock Modelling lockdown measures in epidemic outbreaks using selective
  socio-economic containment with uncertainty.
\newblock {\em Math. Biosci. Eng.}, 18(6):7161--7190, 2021.

\bibitem{Baden2021403}
L.~Baden, H.~El~Sahly, B.~Essink, K.~Kotloff, et~al.
\newblock Efficacy and safety of the m{R}{N}{A}-1273 {S}{A}{R}{S}-{C}o{V}-2
  vaccine.
\newblock {\em New England Journal of Medicine}, 384(5):403--416, 2021.

\bibitem{MR4385929}
E.~Bernardi, L.~Pareschi, G.~Toscani, and M.~Zanella.
\newblock Effects of vaccination efficacy on wealth distribution in kinetic
  epidemic models.
\newblock {\em Entropy}, 24(2):Paper No. 216, 22, 2022.

\bibitem{BertagliaLiuParechi}
G.~Bertaglia, L.~Liu, L.~Pareschi, and X.~Zhu.
\newblock Bi-fidelity stochastic collocation methods for epidemic transport
  models with uncertainties.
\newblock {\em Networks and Heterogeneous Media}, 17(3):401--425, 2022.

\bibitem{MR4009539}
B.~Buonomo, R.~Della~Marca, and A.~d'Onofrio.
\newblock Optimal public health intervention in a behavioural vaccination
  model: the interplay between seasonality, behaviour and latency period.
\newblock {\em Math. Med. Biol.}, 36(3):297--324, 2019.

\bibitem{MauroSIR}
R.~M. Colombo and M.~Garavello.
\newblock Well posedness and control in a nonlocal {S}{I}{R} model.
\newblock {\em Appl Math Optim}, 84:737--771, 2021.

\bibitem{preprint}
R.~M. Colombo, M.~Garavello, F.~Marcellini, and E.~Rossi.
\newblock An age and space structured {S}{I}{R} model describing the
  {C}{O}{V}{I}{D}-19 pandemic.
\newblock {\em Journal of Mathematics in Industry}, 10(1), 2020.

\bibitem{futuro}
R.~M. Colombo, M.~Garavello, F.~Marcellini, and E.~Rossi.
\newblock General renewal equations motivated by biology and epidemiology.
\newblock {\em Preprint}, 2022.

\bibitem{FrancescaElena}
R.~M. Colombo, F.~Marcellini, and E.~Rossi.
\newblock Vaccination strategies through intra--compartmental dynamics.
\newblock {\em Networks and Heterogeneous Media}, 17(3):385--400, 2022.

\bibitem{DiGiamberardinoEtAl}
P.~Di~Giamberardino, R.~Caldarella, and D.~Iacoviello.
\newblock Modeling, analysis and control of {C}{O}{V}{I}{D}-19 in {I}taly:
  Study of scenarios.
\newblock {\em Proceedings of the 18th International Conference on Informatics
  in Control, Automation and Robotics, ICINCO 2021}, pages 677--684, 2021.

\bibitem{MR4147945}
G.~Dimarco, L.~Pareschi, G.~Toscani, and M.~Zanella.
\newblock Wealth distribution under the spread of infectious diseases.
\newblock {\em Phys. Rev. E}, 102(2):022303, 14, 2020.

\bibitem{FabbriGozziZanco}
G.~Fabbri, F.~Gozzi, and G.~Zanco.
\newblock Verification results for age-structured models of economic-epidemics
  dynamics.
\newblock {\em J. Math. Econom.}, 93:102455, 11, 2021.

\bibitem{Giordano2020855}
G.~Giordano, F.~Blanchini, R.~Bruno, P.~Colaneri, et~al.
\newblock Modelling the {C}{O}{V}{I}{D}-19 epidemic and implementation of
  population-wide interventions in {I}taly.
\newblock {\em Nature Medicine}, 26(6):855--860, 2020.

\bibitem{MR3821682}
M.~Groppi and R.~Della~Marca.
\newblock Epidemiological models and vaccinations: from {B}ernoulli to the
  present.
\newblock {\em Mat. Cult. Soc. Riv. Unione Mat. Ital. (I)}, 3(1):45--59, 2018.

\bibitem{Kai2021221}
X.~Kai, T.~Xiao-Yan, L.~Miao, L.~Zhang-Wu, et~al.
\newblock Efficacy and safety of {C}{O}{V}{I}{D}-19 vaccines: A systematic
  review.
\newblock {\em Chinese Journal of Contemporary Pediatrics}, 23(3):221--228,
  2021.

\bibitem{keimerPflug}
A.~Keimer and L.~Pflug.
\newblock Modeling infectious diseases using integro-differential equations:
  Optimal control strategies for policy decisions and applications in
  {C}{O}{V}{I}{D}-19.
\newblock Technical report, Friedrich-Alexander-Universit{\"a}t
  Erlangen-Nuernberg, 2020.

\bibitem{Law2020}
K.~Law, K.~Peariasamy, B.~Gill, S.~Singh, et~al.
\newblock Tracking the early depleting transmission dynamics of
  {C}{O}{V}{I}{D}-19 with a time-varying {S}{I}{R} model.
\newblock {\em Scientific Reports}, 10(1), 2020.

\bibitem{LemonMahmoud2005}
S.~Lemon and A.~Mahmoud.
\newblock The threat of pandemic influenza: Are we ready?
\newblock {\em Biosecurity and Bioterrorism}, 3(1):70--73, 2005.

\bibitem{LeVequeBook2002}
R.~J. LeVeque.
\newblock {\em Finite volume methods for hyperbolic problems}.
\newblock Cambridge Texts in Applied Mathematics. Cambridge University Press,
  Cambridge, 2002.

\bibitem{Lumley2021533}
S.~Lumley, D.~O’Donnell, N.~Stoesser, P.~Matthews, et~al.
\newblock Antibody status and incidence of {S}{A}{R}{S}-{C}o{V}-2 infection in
  health care workers.
\newblock {\em New England Journal of Medicine}, 384(6):533--540, 2021.

\bibitem{McQuadeEtAl}
S.~McQuade, R.~Weightman, N.~Merrill, A.~Yadav, et~al.
\newblock Control of {C}{O}{V}{I}{D}-19 outbreak using an extended {S}{E}{I}{R}
  model.
\newblock {\em Mathematical Models and Methods in Applied Sciences}, 2021.

\bibitem{Merow202027456}
C.~Merow and M.~Urban.
\newblock Seasonality and uncertainty in global {C}{O}{V}{I}{D}-19 growth
  rates.
\newblock {\em Proceedings of the National Academy of Sciences of the United
  States of America}, 117(44):27456--27464, 2020.

\bibitem{Mukhopadhyay202193}
L.~Mukhopadhyay, P.~Yadav, N.~Gupta, S.~Mohandas, et~al.
\newblock Comparison of the immunogenicity \& protective efficacy of various
  {S}{A}{R}{S}-{C}o{V}-2 vaccine candidates in non-human primates.
\newblock {\em Indian Journal of Medical Research}, 153(1):93--114, 2021.

\bibitem{Murray1}
J.~D. Murray.
\newblock {\em Mathematical biology. {I}}, volume~17 of {\em Interdisciplinary
  Applied Mathematics}.
\newblock Springer-Verlag, New York, third edition, 2002.
\newblock An introduction.

\bibitem{MR4334820}
N.~Parolini, L.~Dede', P.~F. Antonietti, G.~Ardenghi, et~al.
\newblock {\tt {SUIHTER}}: a new mathematical model for {C}{O}{V}{I}{D}-19.
  {A}pplication to the analysis of the second epidemic outbreak in {I}taly.
\newblock {\em Proc. A.}, 477(2253):Paper No. 20210027, 21, 2021.

\bibitem{Polack20202603}
F.~Polack, S.~Thomas, N.~Kitchin, J.~Absalon, et~al.
\newblock Safety and efficacy of the {B}{N}{T}162b2 m{R}{N}{A}
  {C}{O}{V}{I}{D}-19 vaccine.
\newblock {\em New England Journal of Medicine}, 383(27):2603--2615, 2020.

\bibitem{MR3881874}
A.~Pugliese and F.~Milner.
\newblock A structured population model with diffusion in structure space.
\newblock {\em J. Math. Biol.}, 77(6-7):2079--2102, 2018.

\bibitem{Randolph2020737}
H.~Randolph and L.~Barreiro.
\newblock Herd immunity: Understanding {C}{O}{V}{I}{D}-19.
\newblock {\em Immunity}, 52(5):737--741, 2020.

\bibitem{RegisEtAl}
S.~Regis, S.~Nuiro, W.~Merat, and A.~Doncescu.
\newblock A data-based approach using a multi-group {S}{I}{R} model with fuzzy
  subsets: Application to the {C}{O}{V}{I}{D}-19 simulation in the islands of
  {G}uadeloupe.
\newblock {\em Biology}, 10(10), 2021.

\bibitem{Ripperger2020925}
T.~Ripperger, J.~Uhrlaub, M.~Watanabe, R.~Wong, et~al.
\newblock Orthogonal {S}{A}{R}{S}-{C}o{V}-2 serological assays enable
  surveillance of low-prevalence communities and reveal durable humoral
  immunity.
\newblock {\em Immunity}, 53(5):925--933.e4, 2020.

\bibitem{Russo20212517}
A.~Russo, A.~Decarli, and M.~Valsecchi.
\newblock Strategy to identify priority groups for {C}{O}{V}{I}{D}-19
  vaccination: A population based cohort study.
\newblock {\em Vaccine}, 39(18):2517--2525, 2021.

\bibitem{Shrestha2021.06.01.21258176}
N.~K. Shrestha, P.~C. Burke, A.~S. Nowacki, P.~Terpeluk, and S.~M. Gordon.
\newblock Necessity of {C}{O}{V}{I}{D}-19 vaccination in previously infected
  individuals.
\newblock {\em medRxiv}, 2021.

\bibitem{VerrelliDellaRosa}
C.~Verrelli and F.~Della~Rossa.
\newblock Two-age-structured {C}{O}{V}{I}{D}-19 epidemic model: Estimation of
  virulence parameters to interpret effects of national and regional feedback
  interventions and vaccination.
\newblock {\em Mathematics}, 9(19), 2021.

\bibitem{Wang20201}
J.~Wang, R.~Jing, X.~Lai, H.~Zhang, et~al.
\newblock Acceptance of {C}{O}{V}{I}{D}-19 vaccination during the
  {C}{O}{V}{I}{D}-19 pandemic in {C}hina.
\newblock {\em Vaccines}, 8(3):1--14, 2020.

\bibitem{Whitley20211774}
R.~Whitley, A.~Babiker, L.~Cooper, S.~Ellenberg, et~al.
\newblock Efficacy of the m{R}{N}{A}-1273 {S}{A}{R}{S}-{C}o{V}-2 vaccine at
  completion of blinded phase.
\newblock {\em New England Journal of Medicine}, 385(19):1774--1785, 2021.

\bibitem{Yang20202708}
C.~Yang and J.~Wang.
\newblock A mathematical model for the novel coronavirus epidemic in {W}uhan,
  {C}hina.
\newblock {\em Mathematical Biosciences and Engineering}, 17(3):2708--2724,
  2020.

\bibitem{Zenilman20211761}
J.~Zenilman, R.~Belshe, K.~Edwards, S.~Self, et~al.
\newblock Safety and efficacy of the {B}{N}{T}162b2 m{R}{N}{A}
  {C}{O}{V}{I}{D}-19 vaccine through 6 months.
\newblock {\em New England Journal of Medicine}, 385(19):1761--1773, 2021.

\end{thebibliography}

}

\end{document}